\documentclass[a4paper,12pt]{article}
\usepackage[active]{srcltx}
\usepackage[ascii]{inputenc}
\usepackage[T1]{fontenc}
\usepackage[ngerman,german]{babel}
\usepackage{amsmath}
\usepackage{amssymb,amsfonts,textcomp}
\usepackage{beton}
\usepackage{euler}
\usepackage{eucal}
\usepackage{color}
\usepackage[top=2cm,bottom=2cm,left=2cm,right=2cm,nohead,includefoot,foot=12pt,footskip=26.144882pt]{geometry}
\usepackage{array}
\usepackage{hhline}
\usepackage{hyperref}
\hypersetup{dvips, colorlinks=true, linkcolor=blue, citecolor=blue, filecolor=blue, urlcolor=blue}
\usepackage[dvips]{graphicx}

\makeatletter
\newcommand\arraybslash{\let\\\@arraycr}
\makeatother
\title{}
\begin{document}
{\centering\bfseries\Large
Geburtstage, W\"urfel, Produkte und Karten
\par}

{\centering\bfseries\large
Einige kombinatorisch{}-stochastische Betrachtungen
\par}

{\centering\small {\bfseries Edgar M.~E.~Wermuth, Technische Hochschule
N\"urnberg}\\ Copyright License: CC
BY-NC-SA 4.0\par}

\bigskip

{\small In diesem Artikel, Ausarbeitung zu einem Vortrag an der TH N"urnberg (7.\,Nov.\,2017),
 werden einige elementare, aber etwas kniffligere Aspekte
und Beispiele aus den im Titel genannten Bereichen vorgestellt und
diskutiert, die ich mir im Laufe der Zeit im Rahmen meiner
Vorlesungst\"atigkeit hier in N\"urnberg zurechtlegte. Zum Teil sind diese
\"Uberlegungen einfach dem Wunsch nach etwas interessanterem
Beispielmaterial und der eige\-nen Neugier 
entsprungen, zum Teil wurden sie
durch Fragen Studierender angeregt.} 

\bigskip

{\centering \bfseries
\large I. Mehrfachgeburtstage \par}

\textbf{I. a) \ }Wir beginnen mit einem gel\"aufigen Beispiel, dem
sogenann\-ten \textit{Geburts\-tags-Paradoxon }(Richard von Mises).
Eigentlich ist nichts paradox daran. Es ist eben nur etwas \"uberraschend,
dass schon bei einer recht kleinen zuf\"allig zu\-sam\-men\-gestellten
Gruppe von Leuten die Wahrscheinlichkeit, dass zwei am sel\-ben Tag
Geburtstag haben, gr\"o{\ss}er als 50\% ausf\"allt. Die Wahrscheinlichkeit,
dass $n$ Per\-so\-nen s\"amtlich an \textit{verschiedenen}
Tagen Geburtstag haben, betr\"agt 
\[
 \boldsymbol{w_n=\frac{365\cdot 364
\,\cdots\,(365-n+1)}{365^n}=
\prod_{k=1}^{n-1}\left(1-\frac{k}{365}\right)}.
\]
Als MAXIMA-Befehl: \verb!w(n):=product(1-k/365,k,1,n-1);! \  MAXIMA
bestimmt die ex\-ak\-ten \textit{ra\-tionalen} Werte. Mit
\ \verb!ev(w(n),float);! \ ergibt sich  $w_{22}\approx
0.5242$, \ $w_{23}\approx 0.4927, \ 
w_{41}\approx 0.0968$, $w_{70}\approx0.00084$. Bei nur 23
Personen ist die Wahrschein\-lich\-keit min\-destens eines
Doppelgeburtstages erstmals {\textgreater}50\%, bei 41 Personen
{\textgreater}90\%, bei 70 \ {\textgreater}99.9\%. Auf den ersten Blick
wirklich erstaunlich.

Wir haben aber bei dieser Rechnung einen \textit{Fehler} gemacht: Der
\textit{gre\-gorianische Kalender}, d.h. der gelegentlich auftretende 29.
Februar, wurde ig\-no\-riert. Das korrigieren wir jetzt und gehen aus von
$\displaystyle n$ rein zuf\"allig zusammen\-gestellten Personen, die alle im
Zeitraum vom 01.01.1901 bis zum 31.12.2016 geboren wurden. In diesem
Zeitraum gab es \textit{alle vier Jahre} ein Schaltjahr mit einem 29.
Februar; auch 2000 war -- anders als 1900 und 2100 -- ein Schaltjahr.

Der repr\"asentative sich wiederholende Zeitraum, aus dem rein zuf\"allig die
Ge\-burts\-tage auszu\-w\"ahlen sind, umfasst nun nicht \ $\displaystyle
365$ Tage, sondern \ $\displaystyle 4\cdot 365+1$. Damit gibt es\\[8pt]
$\displaystyle A^{0}_n=4\cdot 365 \,\cdot \,4\cdot 364 \, \cdots \,4\cdot
(365-n+1)$ \ und \ $\displaystyle A^1_n=n
\,\cdot\,4\cdot365\,\cdot\,4\cdot364\,\cdots\,4\cdot(365-n+2)$\\[8pt]
Auswahlen von $\displaystyle n$ verschiede\-nen Geburtstagen \textit{ohne}
eine bzw. \textit{mit} einer Person, \ die am 29. Februar geboren wurde. So
erhalten wir ({\quotedblbase}G`` stehe f\"ur \textit{gre\-go\-ri\-a\-nisch,
}$\displaystyle 1461=4\cdot365+1$)
\boldmath\[
w^{\mathrm{G}}_n=\frac{A^0_n+A^1_n}{1461^n}
=
\left(1-\frac{3(n-1)}{1461}\right)
\prod_{k=0}^{n-2}\left(1-\frac{4 k+1}{1461}\right),
\]\unboldmath{}%
eine fast genauso einfache Formel wie vorher. Es folgt $
w^{\mathrm{G}}_{22}\approx 0.5247$, $
w^{\mathrm{G}}_{23}\approx 0.4931$, $
w^{\mathrm{G}}_{41}\approx 0.0971$. Die Wahrscheinlichkeiten sind
\textit{nur geringf\"ugig gr\"o{\ss}er}; 23, 41, 70 bleiben die An\-zahlen,
bei denen die Doppel\-ge\-burtstags-Wahrschein\-lich\-keit die 50\%-, die
90\%- bzw. die 99.9\%-Schranke \"ubertrifft. Der Quotient $
w_n^{\mathrm{G}}/w_n$ w\"achst \ mono\-ton, aber nur \textit{sehr langsam}:\\
\ $w_1^{\mathrm{G}}=w_1=1$, \ $\displaystyle
\frac{w_{n+1}^{\mathrm{G}}/w_{n+1}}{w_{n}^{\mathrm{G}}/w_{n}}
= 1 +
\frac{1}{1461}\left(\frac{365\cdot488}{(365-n)(488-n)}-1\right)$, 
$\displaystyle\frac{w_{100}^{\mathrm{G}}}{w_{100}}\approx 1.02158695$,
\ $\displaystyle
\frac{w_{101}^{\mathrm{G}}\,w_{100}}{w_{101}\,w_{100}^{\mathrm{G}}}\approx 1.0005$.

Man kann nat\"urlich noch den Einwand erheben, dass ungepr\"uft alle
Ge\-burts\-tage als gleich wahrscheinlich angenommen werden. Aber
Soziobiologisches klammern wir aus; es soll hier um eine blo{\ss}
\textit{kombinatorische} Er\"orterung gehen.

Wir illustrieren das Geburtstags-Paradoxon anhand eines zweiten Beispiels,
bei dem die Gleichwahrscheinlichkeit der F\"alle wirklich unbezweifelbar
gegeben ist: \ \textit{Lotto{}-Ziehungen}! 

Beim Zahlen-Lotto "`6 aus 49"'  gibt`s bekanntlich $
\binom{49}{6}=13983816$ gleich wahr\-schein\-liche Ziehungsergebnisse. Wie
gro{\ss} ist die Wahrscheinlichkeit, dass im Zeitraum vom 01.1.2017 bis zum
31.12.2060 zweimal \textit{dasselbe Ziehungsergeb\-nis} auftritt?

Wir setzen voraus, dass in diesem Zeitraum jeden Mittwoch und jeden Samstag
eine Ziehung stattfindet. Es sind 44 Jahre, darunter 11 Schaltjahre, also 
\boldmath\[
44 \cdot 365 + 11=16071=7\cdot2295.857142\ldots\text{\ 
Tage.}
\]\unboldmath%
Da der 01.1.2017 ein Sonntag war, ist also der 31.12.2060 ein Freitag, und es
gibt im ins Auge gefassten Zeitraum $\displaystyle 2\cdot2295+1=4591$
Ziehungen. Die Wahr\-schein\-lichkeit, dass diese insgesamt unterschiedlich
ausfallen, ist also
\boldmath\[
\prod_{k=1}^{4590}\left(1 - \frac{k}{13983816}\right)=0.47069\ldots
\]\unboldmath%
Mit Wahrscheinlichkeit 53\% wird mindestens ein Ziehungsergebnis ein zweites
Mal auf\-treten und in Zeitungsmeldungen als Sensation gefeiert werden.

Und die Wahrscheinlichkeit $ w$, dass eine solche
Ziehungs-Koinzidenz nur \textit{genau ein\-mal} auftritt im betrachteten
Zeitraum? F\"ur diese gilt offenbar
\boldmath\[w\!=\!\binom{4591}{2}
\frac{13983816\cdot 13983815 \cdots 13979227}{13983816^{4591}}
\!=\!\binom{4591}{2}\!\prod_{k=1}^{4589}\!\left(1\!-\!\frac{k}{13983816}\right)\!
\approx0.4961.
\]\unboldmath%

\textbf{I. b) \ } Wir diskutieren nun kombinatorisch etwas kniffligere und auch
den Re\-chen\-auf\-wand deutlich hochschraubende Fragen: die
Wahrscheinlichkeit von \textit{Drei\-fach- und Vierfach-Geburtstagen}.

Einfachheitshalber lassen wir den Son\-der\-fall des 29. Februars
unber\"ucksichtigt; zun\"achst jedenfalls.

Wie gro{\ss} ist die Wahrscheinlichkeit, dass von $\displaystyle n$ Personen
mindestens \textit{drei} am selben Tag Geburtstag haben?

Es ist wieder einfacher, erst einmal \textit{Gegen}{}-Anzahlen, also die
An\-zahl der F\"alle mit h\"ochstens zweifach auftretenden gleichen
Geburtstagen zu bestimmen. Das sind alle diejenigen Konstellationen, bei
denen \textit{entweder} alle Geburtstage verschieden sind -- \textit{oder}
aber $\displaystyle k$ Geburtstage doppelt auftreten $\displaystyle (1\leq k
\leq \lfloor n/2\rfloor)$, kom\-bi\-niert mit \ $\displaystyle n-2k$
weiteren Einzelgeburtstagen.

Es gibt $\displaystyle \sum_{k=0}^{\lfloor n/2\rfloor}\binom{n}{2k}
\frac{(2 k)!}{2^k\cdot k!}\prod_{i=0}^{n-k-1}(365-i)$ F\"alle mit
$\displaystyle 0$ bis $\displaystyle \left\lfloor
\frac{n}{2}\right\rfloor$Doppelgeburtstagen.

\textit{Begr\"un\-dung:} 
W\"ahle $\displaystyle 2 k$ der Personen $\displaystyle 1$ bis $\displaystyle
n$ aus, die Doppelgeburtstage haben; dann Paar\-num\-mern $\displaystyle
(1),\ldots,(k)$ zweimal zuordnen, Austausch der Paar-Partner und
Um\-nu\-me\-rierung ergibt \textit{dieselbe} Paar\-ver\-tei\-lung; dann den
Paaren wie den Ein\-zelnen der Reihe nach insgesamt $\displaystyle n-k$
verschie\-de\-ne Geburtstage zuordnen.

Als Wahrscheinlichkeit \textit{mindestens eines} Dreifachgeburtstages folgt
\[
\boxed{w_3(n)=1-\sum_{k=0}^{\lfloor n/2 \rfloor}\binom{n}{2 k}
\prod_{i=1}^k\frac{2 i-1}{365}
\prod_{i=1}^{n-k-1}\left(1-\frac{i}{365}\right)=
1-w_n \sum_{k=0}^{\lfloor n/2 \rfloor}\binom{n}{2 k}
\prod_{i=1}^k\frac{2 i-1}{365-(n-i)}}
\]
Die effizientere zweite Variante der Formel gilt nur f\"ur $
w_n>0$, also $ n\leq 365$. Es gilt aber schon $
w_3(300)\approx \boldsymbol{0.99999999997}$, so dass nur viel kleinere
$ n$ interes\-sie\-ren. Schon bei 194 Personen \"ubertrifft die
Wahrscheinlichkeit erstmals die 99.9\%-Schran\-ke. Eine Tabelle einiger
Werte:\\[6pt]
\begin{tabular*}{\linewidth}{@{\extracolsep\fill}| c || c | c | c | c | c | c | c | c | c | c |}
\hline
$\boldsymbol{n}$ & 46 & 47 & 87 & 88 & 100 & 131& 132 & 168 & 194 & 214\\
\hline
$\boldsymbol{w_3(n)}$ &  0.09996 & 0.1062 & 0.499 & 0.511 & 0.64586 & {\bfseries 0.896} & {\bfseries 0.9014} & {\bfseries 0.9902} & {\bfseries .999078} & {\bfseries .999904}\\ 
\hline
\end{tabular*}

Aus $w_3(250)\approx\boldsymbol{0.99999946}$ folgt, dass unter
2 Millionen Hochschulen mit min\-des\-tens 250 Professorinnen und
Professoren nur eine \textit{nicht} drei Dozenten mit dem\-selben Geburtstag
hat, d.h. h\"ochstwahrscheinlich \textit{keine} auf diesem Planeten.

Nun diskutieren wir den Fall des Auftretens von
\textit{Vierfach}{}-Geburtstagen.

Dazu bestimmen wir, bei insgesamt $ n$ Personen, die Anzahl
aller Geburtstags\-ver\-teilungen mit $ j$ Doppelgeburtstagen
und $ k$ Dreifachgeburtstagen, wobei na\-t\"ur\-lich gelten
muss: $ 0\leq j\leq \lfloor n/2\rfloor, \ 0 \leq k \leq
\lfloor(n-2j)/3 \rfloor$.

Wir w\"ahlen zuerst $ 2 j$ Personen aus f\"ur die Paarbildung
und dann von den ver\-blie\-benen $ n - 2 j$ Personen weitere
$ 3 k$ f\"ur die Bildung von Geburtstags-Tripeln. Das ergibt
$ \displaystyle\binom{n}{2 j}\binom{n-2 j}{3 k}$ verschiedene
Auswahlm\"oglichkeiten. Nun sind noch die dabei m\"oglichen verschiedenen
Paar- und Tripel-Bildungen zu z\"ahlen. Die An\-zahl der Paare ist
$ 1\cdot3\cdot5\cdots(2j-1)$ wie im Fall vorher. Den
$ 3 k$ Personen werden die Tripel-Label $ (1)$ bis
$ (k)$ jeweils dreimal zugeordnet; da eine Permutation
glei\-cher Label untereinander (jeweils $ 3!$ M\"oglichkeiten)
und eine Umnumerierung der Tri\-pel ($ k!$ M\"oglichkeiten) die
Tripel-Gesamtheit nicht \"andern, gibt es $\displaystyle \frac{(3
k)!}{6^k\,k!}$ ver\-schie\-dene Tripel bei gegebener Gruppe von
$ 3 k$ Personen. 

Nun sind noch $ j$ Zweifach-, $ k$ Dreifach- und
$ n-2 j-3 k$ Einfachgeburtstags-Da\-ten, also insgesamt
$ n-j-2k$ verschiedene Geburtstagsdaten den (z.B.
lexiko\-gra\-fisch) aufgelisteten Paaren, Tripeln und Einzelnen zuzuordnen.
Die Anzahl der Geburts\-tags\-verteilungen ohne Vierfachgeburtstage ist
somit
\[ \sum_{j=0}^{\lfloor
n/2\rfloor}\binom{n}{2j}\frac{(2j)!}{2^j\,j!}\sum_{k=0}^{\lfloor
(n-2j)/3\rfloor}\binom{n-2j}{3k}\frac{(3k)!}{6^k\,k!}\prod_{i=0}^{n-j-2k-1}(365-i).
\]
Als Wahrscheinlichkeit mindestens eines Vierfachgeburtstages ergibt sich
\[
w_4(n)=1-\sum_{j=0}^{\lfloor n/2\rfloor}\binom{n}{2j}
\prod_{\iota=1}^j\frac{2\iota-1}{365}
\sum_{k=0}^{\lfloor (n-2j)/3\rfloor}\binom{n-2j}{3k}
\prod_{\kappa=1}^k\frac{(3\kappa-2)(3\kappa-1)}{2\cdot365^2}
\prod_{\lambda=1}^{n-j-2k-1}\left(1-\frac{\lambda}{365}\right)
\]
oder (beachte $\displaystyle \prod_{1\leq i\leq
n-j-2k-1}\bigl(1-\frac{i}{365}\bigr)=\prod_{1\leq i\leq n-1}/
\Bigl(\prod_{n-j\leq i\leq n-1}\cdot\prod_{n-j-2k\leq i \leq n-j-1}\Bigr)$
\ im Falle $\displaystyle n\leq 365$)
\[
\boxed{w_4(n)\!=\!1\!-\!w_n\! \!\sum_{j=0}^{\lfloor
n/2\rfloor}\!\!\binom{n}{2j}\!
\prod_{\iota=1}^j\frac{2\iota-1}{365\!-\!(n\!-\!\iota)}\!\!\!\!
\sum_{k=0}^{\lfloor (n-2j)/3\rfloor}\!\!\binom{n\!-\!2j}{3k}\!\!
\prod_{\kappa=1}^k\frac{(3\kappa-2)(3\kappa-1)}
{2(365\!-\!(n\!-\!j\!-\!\kappa))(365\!-\!(n\!-\!j\!-\!k\!-\!\kappa))}}
\]
Eine Formel, deren Auswertung schon einigen Aufwand erfordert. Mit MAXIMA ist
das kein Problem auf heutigen PCs, wenn auch die Rechendauer bei
gr\"o{\ss}e\-ren $ n$ rasch anw\"achst. Die zweite Fassung,
allerdings nur f\"ur $ n\leq 365$ g\"ultig, ist die deutlich
effizientere. Wertetabelle:\\[6pt]
\begin{tabular*}{\linewidth}{@{\extracolsep\fill}| c || c | c | c | c | c | c | c | c | c | c |}
\hline
$\boldsymbol{n}$ & 40 & 50 & 75 & 100 & 114 & 150 & 187 & 200 & 250 & 260\\
\hline
$\displaystyle \boldsymbol{w_4(n)}$ & 0.0017 & 0.0043 & 0.0212 & 0.0636 & 0.103 & {\bfseries .2646} & {\bfseries .5027} & {\bfseries .591} & {\bfseries .867} & {\bfseries .90215}\\
\hline
\end{tabular*}

Ist die Anzahl der Personen also kaum mehr als halb so gro{\ss} wie das Jahr
Tage hat, ist die Wahr\-schein\-lichkeit von Vierfachgeburtstagen
$ >0.5$; $ w_4(315)\approx\boldsymbol{0.9901}$,
$ w_4(365)\approx\boldsymbol{0.99958}$, und die Be\-rech\-nung
von $ w_4(400)\approx\boldsymbol{0.9999783}$ (erste Formel!)
dauerte schon 189.3 Sekunden.

{\small Gerd Riehl, \textit{Alte Geburtstagsprobleme -- neu gel\"ost}, Math.
Semesterber. (2014), 61:215-231, dis\-ku\-tiert Algorithmen zur Berechnung
von Mehrfachgeburtstags-Wahrscheinlichkeiten. Unsere \"Uberlegungen
entstanden \textit{unabh\"angig} von diesem Artikel und auch von der dort
genannten wei\-teren Literatur zum Thema. Die hier angegebenen Formeln
werden in Riehls Artikel nicht ex\-plizit genannt, und der 29. Fe\-bruar
wird bei Riehl grunds\"atzlich ausgeklammert. In Riehls Ta\-bel\-le, a.a.O.
S.231, bleibt eine ganze Reihe der hier berechneten Wahrscheinlichkeiten und
An\-zah\-len (alle \textbf{fett} hervor\-gehobenen) un\-ge\-nannt, weil nach
Angaben des Autors die Kapazit\"at des ein\-gesetzten Computers
\"uberfordert war. Von mir wurde ein simpler PC aus dem Jahre 2012 be\-nutzt
und MAXIMA, ein in Lisp pro\-gram\-mierter Interpreter, keine
Number-Crunching-Platt\-form.} 

\
\\
Die Erstel\-lung der folgenden beiden Plots
mit MAXIMA ben\"otigte ca. 24 \ bzw. 1270 \ Sekunden.

\includegraphics[width=8.4cm,height=6.5cm]{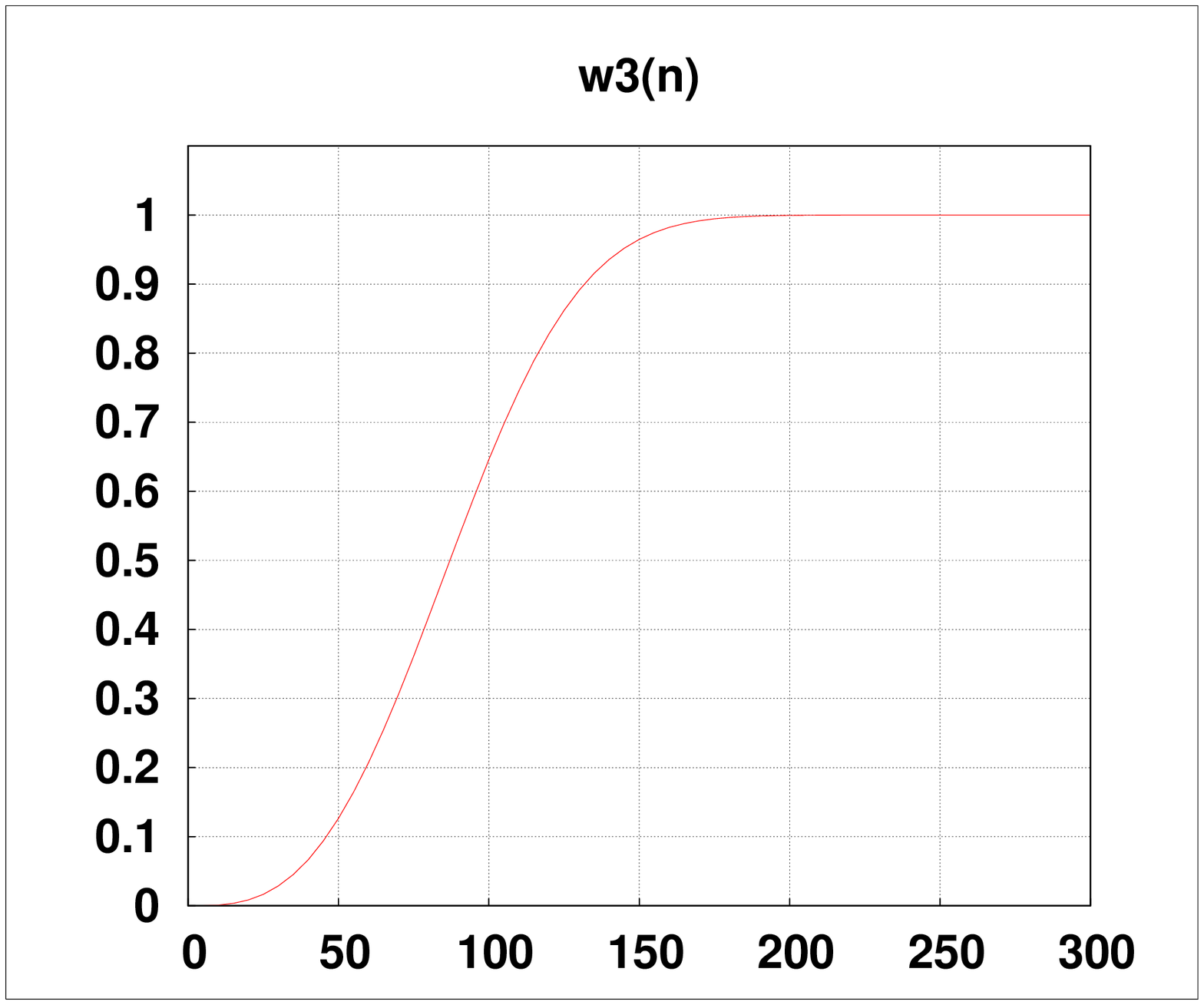}
\ 
\includegraphics[width=8.4cm,height=6.5cm]{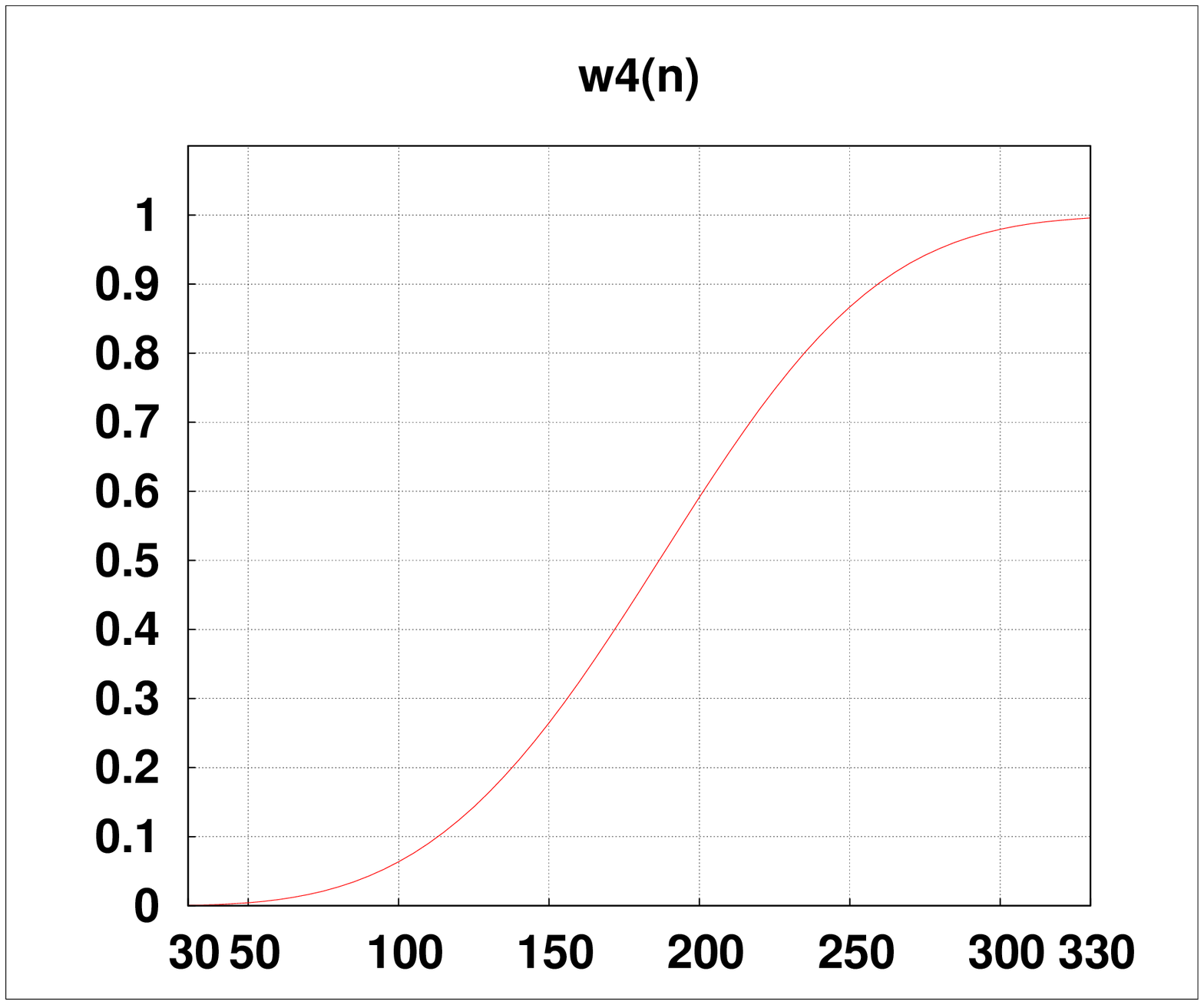}
  
Die verwendeten MAXIMA-Befehle:\\[-24pt]
\begin{verbatim}
w(n):=product(1-k/365,k,1,n-1);    w2(n):=1-w(n);

w3(n):=1-sum(
             binomial(n,2*k)*product((2*i-1)/365,i,1,k)*
             product(1-i/365,i,1,n-k-1),
             k,0,floor(n/2));

w3a(n):=1-w(n)*sum(
                   binomial(n,2*k)*product((2*i-1)/(365-(n-i)),i,1,k),
                   k,0,floor(n/2));
\end{verbatim}
\begin{verbatim}
w4(n):=1-sum(
             binomial(n,2*j)*product((2*i-1)/365,i,1,j)*
             sum(
                 binomial(n-2*j,3*k)*product((3*i-2)*(3*i-1)/(2*365^2),i,1,k)*
                 product(1-i/365,i,1,n-j-2*k-1),
                 k,0,floor((n-2*j)/3)),
             j,0,floor(n/2));
\end{verbatim}
\begin{verbatim}
w4a(n):=1-w(n)*sum(
                   binomial(n,2*j)*product((2*i-1)/(365-(n-i)),i,1,j)*
                   sum(
                       binomial(n-2*j,3*k)*product((3*i-2)*(3*i-1)/
                       (2*(365-(n-j-i))*(365-(n-j-k-i))),i,1,k),
                       k,0,floor((n-2*j)/3)),
                   j,0,floor(n/2));
\end{verbatim}
\begin{verbatim}
load(draw);

draw2d(grid=true,dimensions=[1800,1500],
       point_type=dot, points_joined=true,color=red,
       points(makelist(5*i,i,1,60),makelist(w3a(5*i),i,1,60)),
       yrange=[0.0,1.1],ytics={0.0,0.1,0.2,0.3,0.4,0.5,0.6,0.7,0.8,0.9,1.0},
       xrange=[0,300],font="Helvetica-Bold",font_size=40,
       title="w3(n)",xlabel="",ylabel="",terminal='eps,
       file_name="/home/emew/Dokumente/Geburtstage_u_Co_Juli2017/fig1");
\end{verbatim}

\textbf{I. c) \ } Nun die analogen Formeln f\"ur
\textit{F\"unffach}{}-Geburtstage. Der praktische Re\-chen\-auf\-wand wird
dabei nat\"urlich noch einmal \textit{deutlich} ansteigen.

Wir ermitteln wieder zuerst die \textit{Gegen}{}-Anzahl, also jetzt die aller
Geburtstags\-vertei\-lun\-gen \textit{ohne} F\"unffach-Geburtstage.

Bei insgesamt $\displaystyle n$ Personen geht`s um die Anzahl aller
Geburtstags\-ver\-teilungen mit $\displaystyle j$ Doppelgeburtstagen,
$\displaystyle k$ Dreifachgeburtstagen und $\displaystyle l$
Vierfachgeburtstagen, wobei gelten muss: $\displaystyle 0\leq j\leq \lfloor
n/2\rfloor, \ 0 \leq k \leq \lfloor(n-2j)/3 \rfloor,
0\leq l \leq\lfloor (n-2j-3k)/4\rfloor$.

Es gibt $\displaystyle
\binom{n}{2j}\binom{n-2j}{3k}\binom{n-2j-3k}{4l}$verschiedene
M\"oglichkeiten, bei gegebenen $\displaystyle j$, $\displaystyle k$und
$\displaystyle l$ die Paar-, Tripel- und Quadrupel-Pl\"atze auszuw\"ahlen.
Die Anzahlen dabei m\"oglicher verschiedener Paare und Tri\-pel ken\-nen wir
schon: $\displaystyle \frac{(2j)!}{2^j\,j!}$ bzw. $\displaystyle \frac{(3
k)!}{6^k\,k!}$. Ana\-log ergibt sich als Anzahl verschiedener Quadrupel
$\displaystyle \frac{(4l)!}{24^l\,l!}$.\\ Den Paaren, Tripeln und Quadrupeln
sind $\displaystyle j +k+l+(n-2j-3k-4l)=n-j-2k-3l$ verschiedene Geburtstage
zuzuordnen. Daraus resultieren insgesamt 
\[
\sum_{j=0}^{\lfloor n/2\rfloor}\binom{n}{2j}\frac{(2j)!}{2^j\,j!}
\sum_{k=0}^{\lfloor (n-2j)/3\rfloor}\binom{n-2j}{3k}\frac{(3k)!}{6^k\,k!}
\sum_{l=0}^{\lfloor(n-2j-3k)/4\rfloor}\binom{n-2j-3k}{4l}
\frac{(4l)!}{24^l\,l!}
\prod_{i=0}^{n-j-2k-3l-1}(365-i)
\]
Geburtstagsverteilungen ohne F\"unffach{}-Geburtstage; daher:
{\small\[\begin{split}
\boldsymbol{w_5(n)} =& 1\!-\!\!\sum_{j=0}^{\lfloor
n/2\rfloor}\!\!\binom{n}{2j}\!
\prod_{\iota=1}^j\frac{2\iota\!-\!1}{365}\!
\sum_{k=0}^{\lfloor (n-2j)/3\rfloor}\!\!\binom{n\!-\!2j}{3k}\!
\prod_{\kappa=1}^k\frac{(3\kappa\!-\!2)(3\kappa\!-\!1)}{2\cdot365^2}\\
&\phantom{\quad}\sum_{l=0}^{\lfloor(n-2j-3k)/4\rfloor}\!\!\binom{n\!-\!2j\!-\!3k}{4l}\!
\prod_{\lambda=1}^{l}\frac{(4\lambda\!-\!3)(4\lambda\!-\!2)(4\lambda\!-\!1)}
{6\cdot 365^3}
\prod_{\mu=1}^{n-j-2k-3l-1}\!\!\left(1\!-\!\frac{\mu}{365}\right)
\end{split}
\]}%
Mit 
\[\displaystyle
\boldsymbol{p_1(j)}:=\prod_{\iota=1}^j\frac{2\iota\!-\!1}{365}, \ 
\boldsymbol{p_2(k)}:=\prod_{\kappa=1}^k\frac{(3\kappa\!-\!2)(3\kappa\!-\!1)}
{2\cdot365^2}, \ \boldsymbol{p_3(l)}:=\prod_{\lambda=1}^l
\frac{(4\lambda\!-\!3)(4\lambda\!-\!2)(4\lambda\!-\!1)}{6\cdot365^3}
\]
und $\displaystyle \boldsymbol{p_0(m)}:=\prod_{\mu=1}^m\left(1-\frac{\mu}{365}\right)$
ergibt sich die etwas \"ubersichtlichere Darstellung
{\small\[
\boxed{\boldsymbol{w_5(n)}=1\!-\!\sum_{j=0}^{\lfloor
n/2\rfloor}\!\binom{n}{2j}\,p_1(j)\!
\sum_{k=0}^{\lfloor (n-2j)/3\rfloor}\!\binom{n\!-\!2j}{3k}\,p_2(k)\!
\sum_{l=0}^{\lfloor (n-2j-3k)/4\rfloor}\!\!\binom{n\!-\!2j\!-\!3k}{4l}\,p_3(l)\,
p_0(n\!-\!j\!-\!2k\!-\!3l\!-\!1)}
\]}%
Die Verallgemeinerung auf $m$-fache Geburtstage ist
offensichtlich; dabei kom\-men dann auch $\displaystyle
p_4(k)=\prod_{i=1}^k\frac{(5i\!-\!4)(5i\!-\!3)(5i\!-\!2)(5i\!-\!1)}{24\cdot365^4},
\ p_5(k), \ldots$ ins Spiel.

Der Aufwand bei den Auswertungen der $w_5$-Formel wird etwas
im Zaum ge\-hal\-ten durch \textit{Vorab{}-Berechnung} der Faktoren
$\displaystyle p_1(j) \ (0\leq j\leq n/2)$, $\displaystyle p_2(k) \ (0\leq k
\leq n/3)$, \ $\displaystyle p_3(l) \ (0\leq l\leq n/4)$ sowie $p_0(k)$, die dann
nur noch aufgerufen werden.

{\small Der \textit{einzige} Wert, der im zitierten Artikel von Riehl angegeben ist,
ist das klein\-ste $n$, f\"ur welches $\displaystyle
w_5(n)\!>\!0.01$. Es ist $ w_5(124)\!\approx\!0.009636$,
$\displaystyle w_5(125)\!\approx\!0.010013$.}

Wir stellen einige Rechenergebnisse grafisch dar:

\includegraphics[width=8.4cm,height=7cm]{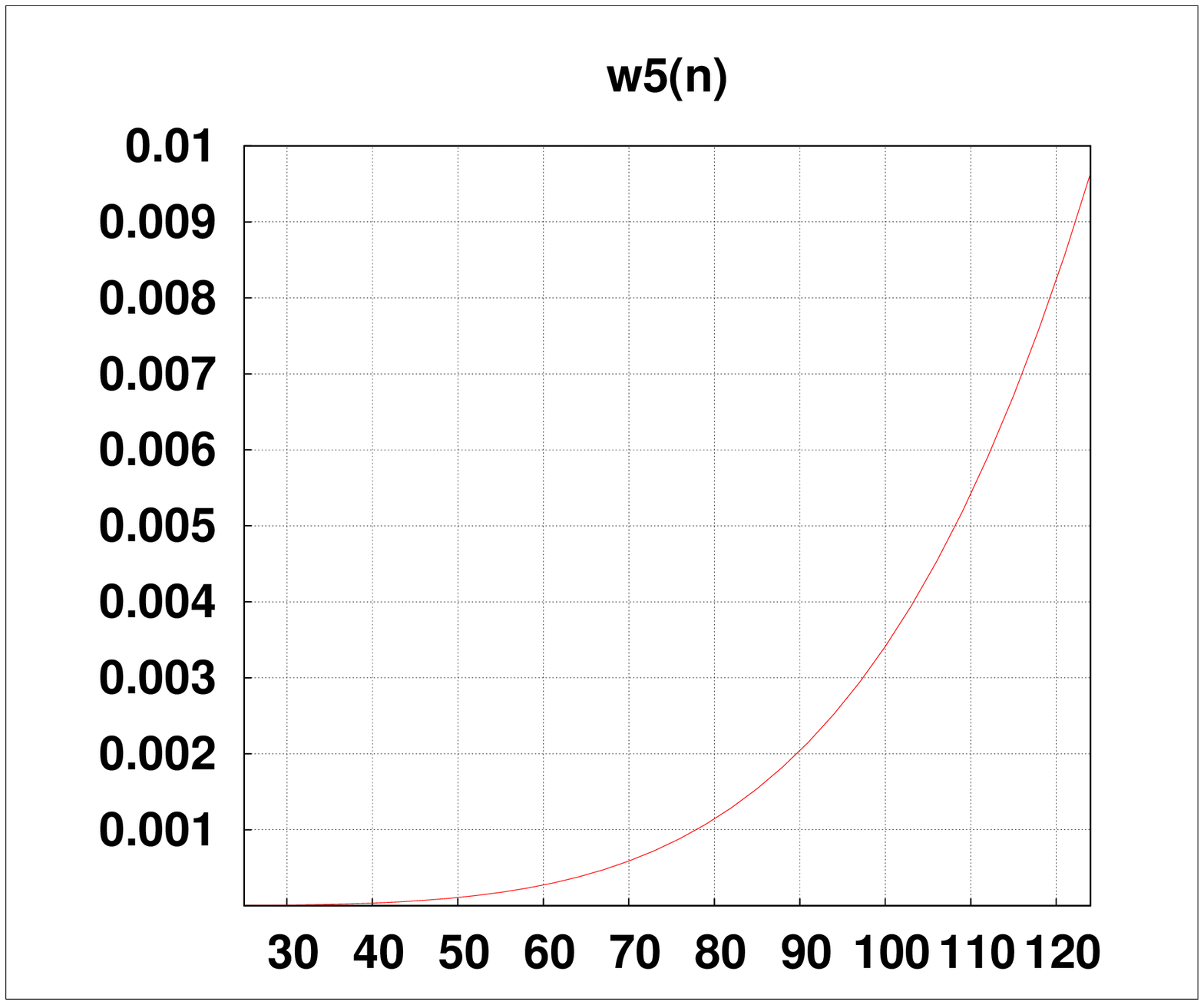}
\ 
\includegraphics[width=8.4cm,height=7cm]{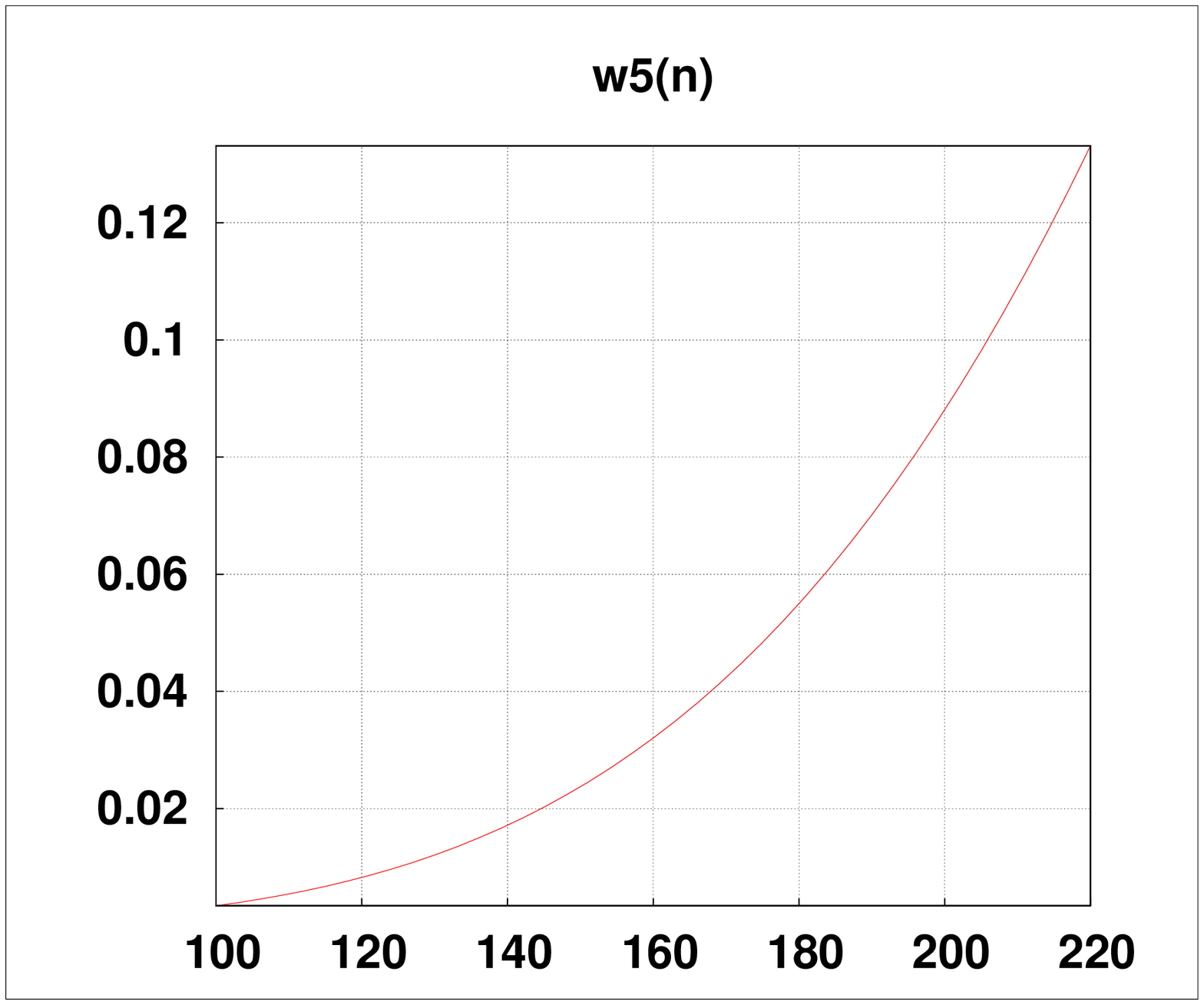}

Die erste MAXIMA-Grafik zeigt die Werte von $ w_5(n)$ f\"ur
$ 25\leq n\leq 124$. Wegen $
w_5(25)\approx0.0000028597435$ sind kleinere $ n$
uninteressant. Die Erstellung der zweiten Grafik mit den 41 $
w_5(n)$-Werten f\"ur $ n=100+3k \, (0\leq k \leq 40)$
er\-forderte ca. 560 Sekunden; $
w_5(206)\approx\boldsymbol{0.100}2585944$. \ \ Das kleinste $
n$ mit $\displaystyle w_5(n)\geq 0.5$ ist $ \boldsymbol{313}$,
da $ w_5(313)\approx\boldsymbol{0.50107}$, $
w_5(312)\approx\boldsymbol{0.496}19557$. Ferner $
w_5(348)\approx\boldsymbol{0.66887}$, $
w_5(347)\approx\boldsymbol{0.664}$.

Der MAXIMA-Code zu obigen Grafiken:\\[-24pt]
{\small\begin{verbatim}
p1(k):=product((2*i-1)/365,i,1,k);   p2(k):=product((3*i-2)*(3*i-1)/(2*365^2),i,1,k);
p3(k):=product((4*i-3)*(4*i-2)*(4*i-1)/(6*365^3),i,1,k); p0(k):=product(1-i/365,i,0,k);

A:makelist(ev(p1(k),bfloat),k,0,300)$   B:makelist(ev(p2(k),bfloat),k,0,200)$
C:makelist(ev(p3(k),bfloat),k,0,150)$   D:makelist(ev(p0(k),bfloat),k,0,600)$

w5(n):=1-sum(binomial(n,2*j)*A[j+1]*
             sum(binomial(n-2*j,3*k)*B[k+1]*
                 sum(binomial(n-2*j-3*k,4*l)*C[l+1]*D[n-j-2*k-3*l],
                     l,0,floor((n-2*j-3*k)/4)),
                 k,0,floor((n-2*j)/3)),}
             j,0,floor(n/2));

load(draw);

draw2d(grid=true,dimensions=[1800,1500],point_type=dot,points_joined=true,color=red,
       points(makelist(100+3*i,i,0,40),makelist(w5(100+3*i),i,0,40)),
       title="w5(n)",xlabel="",ylabel="",terminal='eps,
       font="Helvetica-Bold",font_size=30,
       file_name="/home/emew/Dokumente/Geburtstage_u_Co_Juli2017/fig4");

draw2d(grid=true,dimensions=[1800,1500],point_type=dot,points_joined=true,color=red,
       points(makelist(25+3*i,i,0,33),makelist(w5(25+3*i),i,0,33)),yrange=[0,0.01],
       ytics={0.001,0.002,0.003,0.004,0.005,0.006,0.007,0.008,0.009,0.01},
       title="w5(n)",xlabel="",ylabel="",terminal='eps,
       font="Helvetica-Bold",font_size=30,
       file_name="/home/emew/Dokumente/Geburtstage_u_Co_Juli2017/fig3");
\end{verbatim}}\vspace*{-0.5cm}

Mit doch schon betr\"achtlichem Zeitaufwand k\"onnen wir auch den weiteren
Ver\-lauf von $w_5$ berechnen und grafisch darstellen:

\includegraphics[width=8.4cm,height=7.0cm]{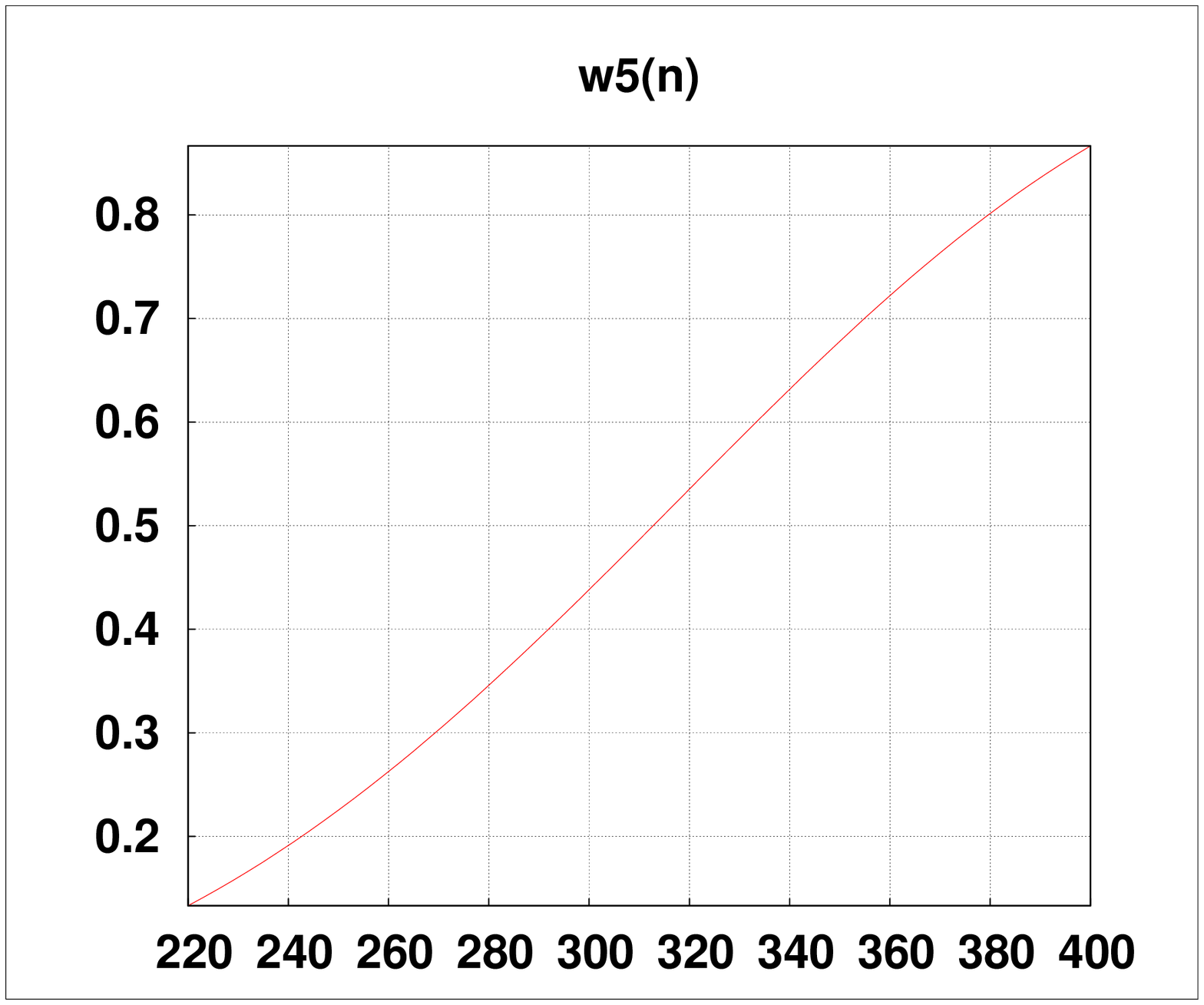}
\   
\includegraphics[width=8.4cm,height=7.0cm]{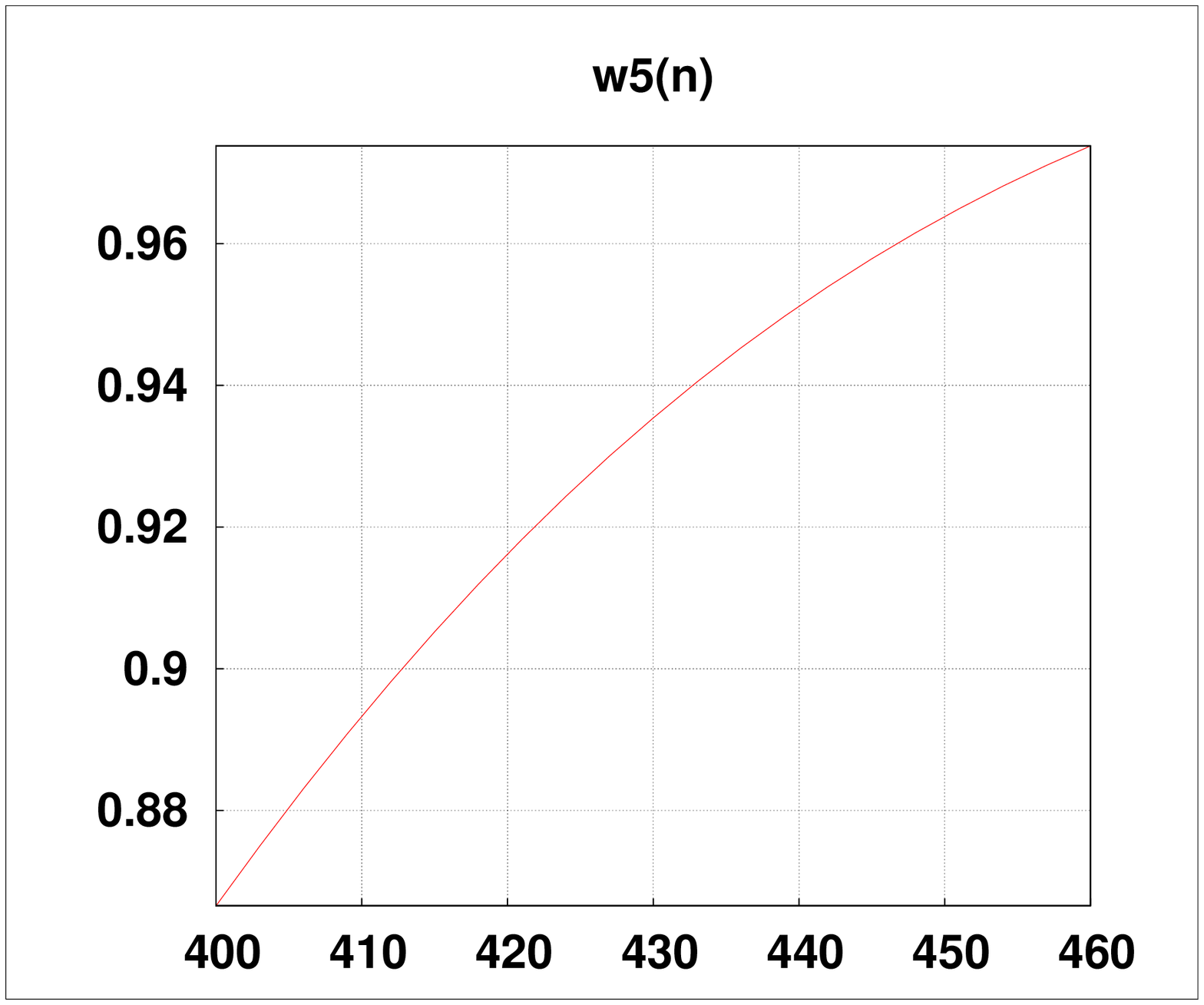}

F\"ur's erste Bild ($w_5(220\!+\!3k)\,(0\!\leq\! k\!\leq\!
60)$) ben\"otigte MAXIMA ca. 7570 Se\-kun\-den, f\"ur`s zweite 7535
($w_5(400\!+\!\!3i)\,(0\!\leq\! i\!\leq\! 20)$) und f\"ur`s
dritte ($\displaystyle w_5(460\!+\!3i)\,(0\!\leq\!i\!\leq\!20)$)
rund 12135.

\includegraphics[width=8.4cm,height=7.5cm]{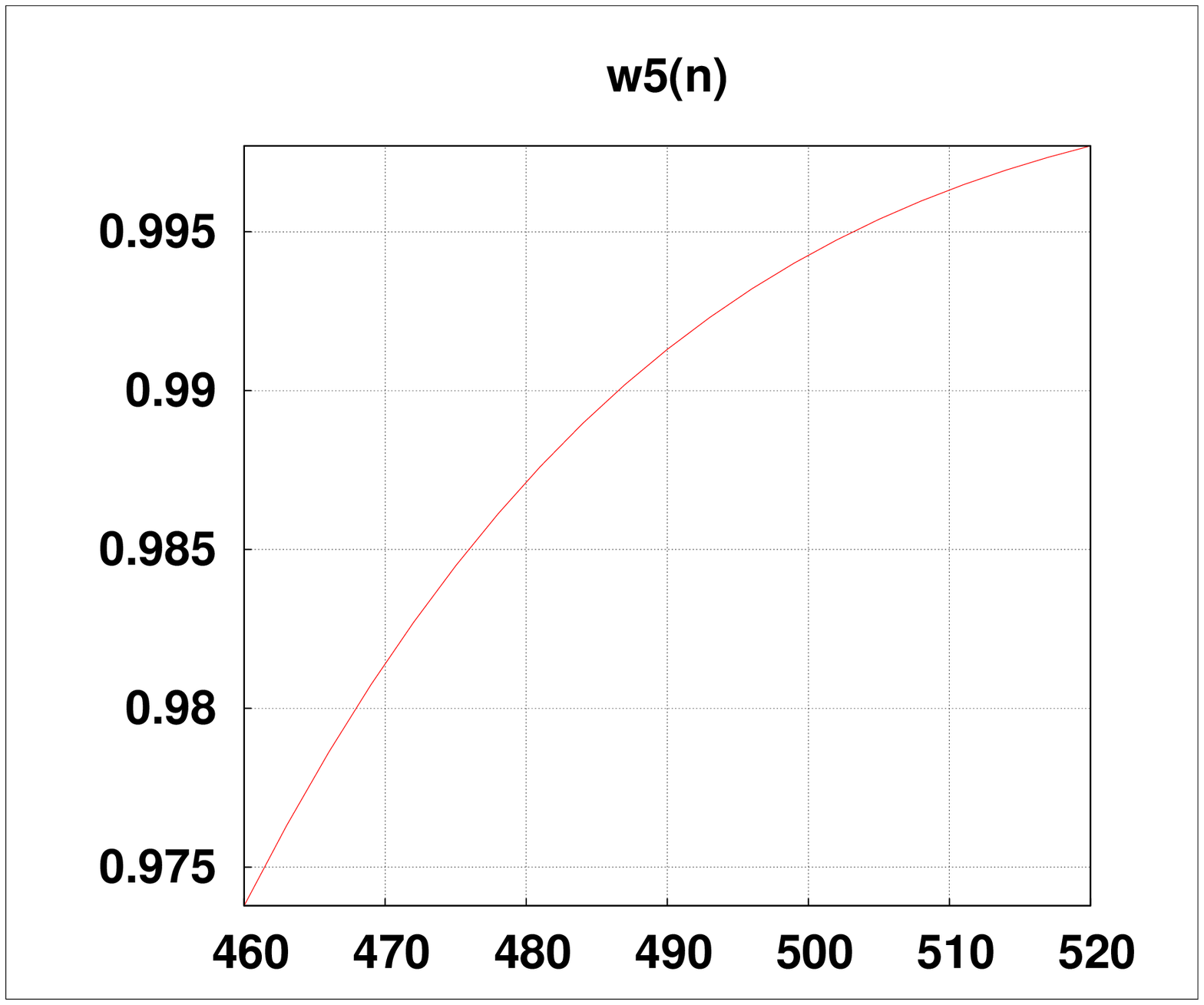}\,
\begin{minipage}[b]{0.49\linewidth}
$\displaystyle w_5(413)\approx\boldsymbol{0.9005775}$,
\  $\displaystyle w_5(412)\approx\boldsymbol{0.89819}$, \ ferner
 \ $\displaystyle w_5(486)\approx\boldsymbol{0.98979}$ \ \ und
$w_5(487)\approx\boldsymbol{0.9901886}$;
 schlie{\ss}lich  $w_5(536)\approx\boldsymbol{.998978}$, $w_5(537)\approx\boldsymbol{0.99903}$. 
Schon $ w_5(487)$ und $\displaystyle w_5(537)$ allein
be\-n\"o\-tigen rund 565 bzw. 810 Sekun\-den an Rechenzeit.\\
{\small Angesichts des doch sehr gro{\ss}en Re\-chen\-zeit-Bedarfs (Verwendung
ein\-fa\-cher \textit{Floats} statt \textit{Big Floats} erg\"abe we\-ni\-ger
als 15\% Verk\"urzung) w\"are es viel\-leicht bes\-ser, bei den
$\displaystyle w_5${}-Berech\-nun\-gen mit ei\-nem echten Numerik-System wie
z.B. SCILAB zu arbei\-ten. \ \ \ \ (Benutzte Software: MAXIMA 5.38.1, Linux
OpenSuse 42.2; Hardware: einfacher Standard-PC, Baujahr 2012.)}
\end{minipage}

Eine zusammenfassende Anzahl-Tabelle zu den $\displaystyle w_k(n)$:

\begin{tabular*}{\linewidth}{@{\extracolsep\fill} |c|cccccc|}
\hline
%\multicolumn{1}{| c |}{ }&
\multicolumn{7}{| c |}{$ \boldsymbol{n_{k,\alpha}:=\min
\{n\in\mathbb{N} \mid w_k(n)\geq \alpha\}}$ }\\
\hline
$\boldsymbol{\alpha}$ & \textbf{0.01} & \textbf{0.1} & \textbf{0.5} & \textbf{0.9} & 
\textbf{0.99} & \textbf{0.999} \\
\hline
$\boldsymbol{k=2}$ &  4 & 10 & 23 & 41 & 57 & 70 \\ 
\hline
$\boldsymbol{k=3}$ &  22 & 47 & 88 & 132 & 168 & 194 \\ 
\hline
$\boldsymbol{k=4}$ &  62 & 114 & 187 & 260 & 315 & 353 \\ 
\hline
$\boldsymbol{k=5}$ &  125 & 206 & 313 & 413 & 487 & 537 \\ 
\hline
\end{tabular*}

{\small (Alle Anzahlen $\displaystyle >125$, also zwei Drittel der nichttrivialen
Werte ($\displaystyle k>2$), blieben im zitierten Artikel von G. Riehl
offen.)}

\bigskip

\textbf{I. d) \ }  Die generelle Ber\"ucksichtigung des \ \textit{29. Februars}.

Zun\"achst noch einmal der Fall von \textit{Doppel}geburtstagen. Mit der
Wahrschein\-lich\-keit  $ p$ einer Geburt am 29.02. gilt \ 
\boldmath\(
w_n^{\mathrm{G}}=(1-p)^nw_n+np(1-p)^{n-1}w_{n-1}. 
\)\unboldmath%

F\"ur $ p=1/1461$ stimmt diese Formel \textit{exakt} mit der
fr\"uheren \"uberein. Man kann durch an\-de\-re $p$-Wahl
sich leicht an eine Auswahl z.B. von zwischen 1880 und 1910 Geborenen
anpassen. (Wichtig ist, das Zufalls\-experi\-ment \textit{pr\"a\-zise zu
defi\-nie\-ren}: 
"`eine rein zuf\"allige Auswahl von Personen,
die im Zeitraum von $\displaystyle t_1$ bis $\displaystyle t_2$ geboren
sind..."'.) \ \ \ F\"ur das Verh\"altnis zu $\displaystyle w_n$ folgt 

{\small$\displaystyle
w_n^{\mathrm{G}}/w_n=(1-p)^{n-1}\left(1+\frac{p(n-1)}{1-n/366}\right)>
(1-(n-1)p)\left(1+\frac{p(n-1)}{1-n/366}\right)=1+\frac{(n-1)p(n/366-(n-1)p)}{1-n/366}$.}

Nun zum Allgemeinfall. Folgende Ereignisse spielen eine Rolle:

$\boldsymbol{F_k}$ : "`$\displaystyle k$ von
$n$ Personen haben am 29. Februar Geburtstag"';

$\boldsymbol{\leq k}$ : "`h\"ochstens
$k$-fache Geburtstage bei einer Gruppe von $\displaystyle n$
Personen"';

$\boldsymbol{\geq k}$ : "`mindestens ein
$k$-facher Geburtstag bei einer Gruppe von $\displaystyle n$
Personen"'.

Dann gilt offenbar (wobei $\displaystyle p
=\frac{0.2425}{365.2425}=\frac{97}{400\cdot 365+97}$ bei
\textit{gregorianischem} Auswahlzeitraum von $400$ Jahren und $p=1/1461$ bei \textit{julianischer} Schaltung
im Auswahlzeitraum)
\boldmath\[ 
\mathrm{P}(F_k)=\binom{n}{k}p^k(1-p)^{n-k}\ (0\leq
k\leq n), \ \ \mathrm{P}(\leq k \mid
F_i)=1-w_{k+1}(n-i) \ (i\leq k).
\]
Beispielsweise die "`gregorianische"' Formel f\"ur die
Wahrscheinlichkeit min\-des\-tens eines \textit{Drei\-fach}{}-Geburtstags
ist (Formel von der totalen Wahrscheinlichkeit)
\[
\mathrm{P}(\geq 3)=1 -
 \mathrm{P}(\leq 2 \mid F_0)\cdot\mathrm{P}(F_0)-
\mathrm{P}(\leq 2 \mid F_1)\cdot\mathrm{P}(F_1)-
\mathrm{P}(\leq 2 \mid F_2)\cdot\mathrm{P}(F_2).
\]\unboldmath%
und daher 
\[
\boxed{\!\boldsymbol{w_3^{\mathrm{G}}(n)}\!=\!1\!-\!(1\!-\!p)^n\bigl(1\!-\!w_3(n)\bigr)
\!-\!n p (1\!-\!p)^{n-1}\bigl(1\!-\!w_3(n\!-\!1)\bigr)
\!-\!{\scriptstyle\frac{n(n-1)}{2}}p^2(1\!-\!p)^{n-2}
\bigl(1\!-\!w_3(n\!-\!2)\bigr)\!\!}
\]
V\"ollig analog ergeben sich in allen anderen F\"al\-len mehr\-fa\-cher
Geburtstags-Ko\-in\-zidenzen die gregorianischen Formeln aus den
schalttagfreien.

\begin{tabular*}{\linewidth}{@{\extracolsep\fill} |c|cccccc|}
\hline
%\multicolumn{1}{| c |}{ }&
\multicolumn{7}{| c |}{$ \boldsymbol{n^G_{k,\alpha}:=\min
\{n\in\mathbb{N} \mid w^G_k(n)\geq \alpha\}}$ }\\
\hline
$\boldsymbol{\alpha}$ & \textbf{0.01} & \textbf{0.1} & \textbf{0.5} & \textbf{0.9} & 
\textbf{0.99} & \textbf{0.999} \\
\hline
$\boldsymbol{k=2}$ &  4 & 10 & 23 & 41 & 57 & 70 \\ 
\hline
$\boldsymbol{k=3}$ &  22 & 47 & 88 & 132 & 168 & 194 \\ 
\hline
$\boldsymbol{k=4}$ &  62 & 114 & 187 & 260 & 315 & 354 \\ 
\hline
$\boldsymbol{k=5}$ &  126 & 206 & 313 & 413 & 487 & 537 \\ 
\hline
\end{tabular*}

{\small Dabei $\boldsymbol{p=1/1461}$. Nur \textit{zwei} Tabellenwerte unterscheiden sich von den entsprechenden
schalttagfreien, sind um 1 gr\"o{\ss}er:
$w_5^{\mathrm{G}}(125)\approx0.009980994$, \ $w_5^{\mathrm{G}}(126)\approx0.010368$, \ 
$w_5(125)\approx0.010013$; \ \ $w_4^{\mathrm{G}}(353)\approx0.99899$, \ 
$w_4^{\mathrm{G}}(354)\approx0.999058$, \ $w_4(353)\approx0.999007595$.}

\medskip

Noch ein paar naheliegende Anmerkungen zur \textit{Monotonie} der
Wahrscheinlich\-kei\-ten $w_k(n)$ und $w_k^{\mathrm{G}}(n)$.

Bei den Gegenwahrscheinlichkeiten $w_n=1-w_2(n)$ und
$w_n^{\mathrm{G}}=1-w_2^{\mathrm{G}}(n)$ ist unmittelbar die
fallende Monotonie abzulesen, und auch $w_2(n)\!>\!w_2^G(n)$
ist klar, da offenbar \ $w_n^G/w_n>1$ im relevanten Bereich
$n\leq 365$.

Die angegebenen For\-meln f\"ur $w_3(n), w_3^{\mathrm{G}}(n),
w_4(n), \ldots$ hingegen sind zu ver\-wickelt, um die Monotonie unmittelbar
abzulesen.

Wir k\"onnen aber die Monotonie ganz einfach \textit{kombinatorisch}
erschlie{\ss}en:

Sei $A_k(n)$ die Anzahl der Geburtstagsverteilungen f\"ur
$n$ Personen mit min\-des\-tens einem $k${}-fachen Geburtstag (sogenannte $k$-\textit{Verteilungen}). 
Ordnen wir einer $(n+1)$-ten Person irgendeinen Geburtstag zu, erhalten wir aus jeder $k$-Vertei\-lung f\"ur $n$ Personen eine $k$-Verteilung f\"ur $n+1$ Personen.

Es folgt also offensichtlich $A_k(n+1)>365\cdot A_k(n)$, da ja
auch aus manchen Nicht-$k$-Verteilungen durch die Hinzunahme
des $(n+1)$-ten Geburtstags $k$-Ver\-tei\-lungen werden. Also \ \ 
\boldmath$\displaystyle w_k(n+1)=\frac{A_k(n+1)}{365^{n+1}}>w_k(n)$\unboldmath.

Dieses Argument, mit 1461 anstelle von 365, funktioniert auch im Falle
"`julia\-nisch"' korrigierter Wahrscheinlichkeitsberechnungen. 

Hat man einen \textit{beliebigen} Zeitraum, aus dem irgendwelche
$n$ gleichwahr\-schein\-lichen Geburtstage ausgew\"ahlt sind,
hat nicht nur der 29. Februar eine ab\-weichende Wahrscheinlichkeit. 

Z.B. bei einem Zeitraum, der am 01.3. eines Jah\-res beginnt und am 30.7.
etliche Jahre sp\"ater endet, sind Geburtstage vom 01.3. bis 30.7.
geringf\"ugig wahrscheinlicher als Geburtstage im Januar und Fe\-bruar bzw.
im August bis Dezember. 

Davon ausgehend, dass jeder Tag des Ge\-samtzeitraumes gleichwahrscheinlich
als Geburtstag infragekommt, kann man mit der \textit{Gesamtanzahl} der Tage
an\-stelle von 365 genauso kombinatorisch argu\-mentieren. Denn
trivialerweise gilt weiterhin: Jede $k$-Verteilung von
$n$ Geburts\-tagen ergibt durch Hinzuf\"ugen
\textit{irgend\-eines} weiteren Geburtstages eine $k$-Ver\-teilung von $n+1$ Geburtstagen. Also gilt ganz
allgemein
\boldmath$\displaystyle w_k^G(n+1) > w_k^G(n)$\unboldmath.

\clearpage

{\centering\bfseries\large II. W\"urfel-Wahrscheinlichkeiten\par}

\textbf{II. a) \ } Zuerst befassen wir uns mit den Wahrscheinlichkeiten f\"ur die
einzelnen er\-ziel\-ba\-ren Augenzahlen beim \textit{W\"urfeln mit mehreren
gleichen W\"urfeln}.

Im Falle \textit{zwei\-er} W\"urfel ist es ein kombinatorisch trivialer
Sachverhalt: 

Ist $ X_2\!\in\!\{2, 3, \ldots, 12\}$ die zugeh\"orige
Augenzahl, gilt \ $\displaystyle \mathrm{P}(X_2\!=\!k)=\frac{6-|7\!-\!k|}{36}
\ (2\leq k\leq 12)$.

Aber schon bei \textit{drei} W\"urfeln werden die Einzelheiten etwas
verwickelter: 

Auf wie viele Weisen beispielsweise kommt die Augenzahl 9 zustande? 

{\centering
$\displaystyle 6+2+1$, \ \ $\displaystyle 5+3+1$, \ \ $\displaystyle 5+2+2$,
\ \ $\displaystyle 4+4+1$, \ \ $\displaystyle 4+3+2$, \ \ $\displaystyle
3+3+3$
\par}

sind die verschiedenen Summen-Muster, denen, der Reihe nach, $
6, \ 6, \ 3, \ 3, \ 6$ und $1$ Elemen\-tar\-ergebnisse beim
Ausf\"uhren dieses W\"urfelexperi\-ments entspre\-chen. Die
Ge\-samt\-wahrscheinlich\-keit der Augenzahl 9 ist somit \ $\displaystyle
\frac{25}{216}$.

Die Anzahl der verschiedenen Muster ist leicht mit der
\textit{Multimengen}-Formel zu erfassen: Wir haben 6 Objekte (die
Augenzahlen) und bilden dreielementige Men\-gen (die W\"urfe), bei denen
auch Objekte mehrfach vorkommen d\"urfen. Die Anzahl der Augenzahl-Muster
ist also $\displaystyle \binom{6\!+\!3\!-\!1}{3}=56$. Die Anzahl der
verschie\-de\-nen und gleichwahrscheinlichen W\"urfe hingegen ist
$6^3=216$, weil alle drei W\"ur\-fel als unabh\"angig und
gleichartig funktionierende Zufallszahlengenerator-In\-di\-vi\-duen zu
betrachten sind. (Die "`Teilchen"'  der Quantenmechanik haben
keinen solchen eindeutigen Individualcharakter, weshalb es die
Bose/Einstein- und Fermi/Dirac-Statistik gibt.)

F\"ur alle Augenzahlen $k$ von 3 bis 18 einzeln festzustellen,
wie viele verschie\-de\-ne W\"urfe jedem der zugeh\"origen Muster
entsprechen, um so die jeweilige Ge\-samt\-zahl $
\boldsymbol{N_3(k)}:=\#(X_3=k)$ zu ermitteln, w\"are zu umst\"andlich.
Deshalb jetzt sys\-te\-ma\-ti\-sche\-re \"Uberlegungen zur Bestimmung aller
$ N_3(k)$. Es gilt offenbar:
\[\begin{split}
&\sum_{k=3}^{18} N_3(k)\cdot x^k 
=(x+x^2+x^3+x^4+x^5+x^6)^3=
x^3\!+\!3x^4\!+\!6x^5\!+\!10x^6\!+\!15x^7\!+\!21x^8\!+\!25x^9\\
&\phantom{\sum_{k=3}^{18} N_3(k)\cdot x^k 
={}}+\!27x^{10}
\!+\!27x^{11}\!
 +\!25x^{12}\!+\!21x^{13}\!+\!15x^{14}\!+\!10x^{15}\!+\!6x^{16}
\!+\!3x^{17}\!+\!x^{18}
\end{split}\]
Will man nicht ein CA-System wie MAXIMA f\"ur's Ausmultiplizieren nutzen,
kann man die $N_3(k)$ auch nach Art des Pascalschen Dreiecks
\textit{von Hand} berechnen. Dazu folgende Umformung:
\[
(x+x^2+x^3+x^4+x^5+x^6)^3 
=\Bigl(x^2(1+x^3)\bigl(\frac{1}{x}+1+x\bigr)\Bigr)^3
= (x^6 +3x^9+3x^{12}+x^{15})\Bigl(\frac{1}{x}+1+x\Bigr)^3
\]
Die letzte Faktorisierung legt die Berechnung der $N_3(k)$ in
drei Schritten nahe. Wegen des dreifachen Faktors $(1/x+1+x)$
wird in der zweiten bis vierten Zeile immer die Summe der drei \"uber einer
Stelle stehenden Zahlen gebildet: $\displaystyle
\searrow\,\downarrow\,\swarrow$
\begin{center}
\begin{tabular}{|>{$}c<{$}|>{$}c<{$}|>{$}c<{$}|
>{$}c<{$}|>{$}c<{$}|>{$}c<{$}|>{$}c<{$}|>{$}c<{$}|>{$}c<{$}|>{$}c<{$}|
>{$}c<{$}|>{$}c<{$}|>{$}c<{$}|>{$}c<{$}|>{$}c<{$}|>{$}c<{$}|}
\hline
3&4&5&6&7&8&9&10&11&12&13&14&15&16&17&18\\
\hline\hline
&&&1&&&3&&&3&&&1&&&\\
&&1&1&1&3&3&3&3&3&3&1&1&1&&\\
&1&2&3&5&7&9&9&9&9&7&5&3&2&1&\\
\boldsymbol{1}&\boldsymbol{3}&\boldsymbol{6}&\boldsymbol{10}&\boldsymbol{15}&\boldsymbol{21}
&\boldsymbol{25}&\boldsymbol{27}&\boldsymbol{27}&\boldsymbol{25}&\boldsymbol{21}&\boldsymbol{15}
&\boldsymbol{10}&\boldsymbol{6}&\boldsymbol{3}&\boldsymbol{1}\\
\hline
\end{tabular}
\end{center} 
Die Augenzahlen $\displaystyle 10$ und $\displaystyle 11$ machen zusammen
\textit{genau ein Viertel} aller F\"alle aus, $7$ \ und
$8$ bzw. $13$ und $14$ zusammen jeweils \textit{ein Sechstel}.

Nun \textit{vier} W\"urfel. Die Anzahl der verschiedenen
gleichwahrscheinlichen W\"urfe ist $6^4=1296$, die Anzahl der
verschiedenen Augenzahl-Muster $\displaystyle \binom{6\!+\!4\!-\!1}{4}=126$.
F\"ur die Berechnung der $\displaystyle N_4(k)\ (4\leq k\leq 24)$ von Hand
ist die Umformung
\[
(x+x^2+x^3+x^4+x^5+x^6)^4=(x^8+4x^{11}+6x^{14}+4x^{17}+x^{20})\Bigl(\frac{1}{x}+1+x\Bigr)^4
\]
ma{\ss}gebend. Das resultierende Schema zeigt u.a. $
\mathrm{P}(10\!\leq\! k\! \leq\! 18)=\frac{29}{36}\approx0.80555...$ :
{\small\begin{center}
\begin{tabular}{|@{\,\,\,}>{$}c<{$}@{\,\,\,}|@{\,\,\,}>{$}c<{$}@{\,\,\,}|@{\,\,\,}
>{$}c<{$}@{\,\,\,}|@{\,\,\,}
>{$}c<{$}@{\,\,\,}|@{\,\,\,}>{$}c<{$}@{\,\,\,}|@{\,\,}>{$}c<{$}@{\,\,}|@{\,\,}>{$}c<{$}@{\,\,}|@{\,\,}
>{$}c<{$}@{\,\,}|@{\,\,}>{$}c<{$}@{\,\,}|@{\,\,}>{$}c<{$}@{\,\,}|@{\,\,}
>{$}c<{$}@{\,\,}|@{\,\,}>{$}c<{$}@{\,\,}|@{\,\,}>{$}c<{$}@{\,\,}|@{\,\,}>{$}c<{$}@{\,\,}|@{\,\,}
>{$}c<{$}@{\,\,}|@{\,\,}>{$}c<{$}@{\,\,}|@{\,\,}>{$}c<{$}@{\,\,}|@{\,\,}>{$}c<{$}@{\,\,}|@{\,\,}
>{$}c<{$}@{\,\,}|@{\,\,}>{$}c<{$}@{\,\,}|@{\,\,}>{$}c<{$}@{\,\,}|}
\hline
4&5&6&7&8&9&10&11&12&13&14&15&16&17&18&19&20&21&22&23&24\\
\hline\hline
&&&&1&&&4&&&6&&&4&&&1&&&&\\
&&&1&1&1&4&4&4&6&6&6&4&4&4&1&1&1&&&\\
&&1&2&3&6&9&12&14&16&18&16&14&12&9&6&3&2&1&&\\
&1&3&6&11&18&27&35&42&48&50&48&42&35&27&18&11&6&3&1&\\
\boldsymbol{1}&\boldsymbol{4}&\boldsymbol{10}&
\boldsymbol{20}&\boldsymbol{35}&\boldsymbol{56}&
\boldsymbol{80}&\boldsymbol{104}&\boldsymbol{125}&
\boldsymbol{140}&\boldsymbol{146}&\boldsymbol{140}&
\boldsymbol{125}&\boldsymbol{104}&\boldsymbol{80}&
\boldsymbol{56}&\boldsymbol{35}&\boldsymbol{20}&
\boldsymbol{10}&\boldsymbol{4}&\boldsymbol{1}\\
\hline
\end{tabular}
\end{center} }
Wir k\"onnen mit der Potenzreihen-Methode auch eine \textit{allgemeine}
Formel auf\-stel\-len f\"ur $\displaystyle N_n(k)$, die Anzahl der
$\displaystyle n${}-W\"urfel-W\"urfe mit Augensumme $\displaystyle k$:
\[\begin{split}
&(x + x^2 + \cdots + x^6)^n =x^n(1+ \cdots +x^5)^n=x^n (1-x^6)^n(1-x)^{-n}\\
&{} = x^n \cdot \sum_{k=0}^n \binom{n}{k}(-1)^kx^{6k} 
\cdot \sum_{l=0}^{\infty} (-1)^l\binom{-n}{l}x^l
= x^n \cdot \sum_{m=0}^{5n} x^m \cdot \sum_{\substack{0\leq k \leq n,\,\, l
\geq 0\\ 6k+l=m}}
\binom{n}{k} \binom{-n}{l}(-1)^{k+l}\\
&{} = \sum_{m=0}^{5n} x^{n+m} \cdot \sum_{\substack{0\leq k \leq \lfloor
m/6\rfloor\\l=m-6k }} \binom{n}{k} 
\binom{-n}{l} (-1)^{k+l}
\end{split}
\]
Bei der vorletzten Zweifachsumme (Summationsgrenze $5n$) wird
ber\"ucksichtigt, dass nur Potenzen bis $x^{6n}$ einen von
$0$ verschiedenen Koeffizienten haben k\"on\-nen. Es ergibt
sich die f\"ur maschinelle Berechnung taugliche Formel:
\boldmath\[
\boxed{
 N_n(k) = \sum_{i=0}^{\lfloor (k-n)/6\rfloor} \binom{n}{i} \binom{-n}{k-n-6i}
\cdot (-1)^{k-n+i}
\quad( n\leq k \leq 6n)}
\]\unboldmath%

\bigskip

\textbf{II. b) \ } Nun w\"urfeln wir mit nur \textit{einem} W\"urfel. Und zwar
bis zum ersten Auf\-tre\-ten einer Sechs, bis zur ersten Doppelsechs (= zwei
Sechsen unmittelbar nach\-ein\-ander), bis zur ersten $\displaystyle
n${}-fach-Sechs. Dem entsprechen die Zufallsgr\"o{\ss}en $\displaystyle
X_1$, $\displaystyle X_2$, $\displaystyle X_n$, welche die Anzahl der
jeweils n\"otigen W\"urfe angeben.

Wie lange muss man im Schnitt w\"urfeln, um eine Doppelsechs zu erleben? Wie
stark streut diese Zahl? Es geht also um \textit{Erwartungswerte} und
\textit{Varianzen} dieser Zufallsgr\"o{\ss}en.

Bei $X_1$ ist es noch simpel: 

$\displaystyle
\mathrm{P}(X_1\!=\!k)=\!\Bigl(\frac{5}{6}\Bigr)^{k-1}\frac{1}{6}$, also
$\displaystyle \mathrm{E}(X_1)=\!
\frac{1}{6}\sum_{k=1}^{\infty} k \Bigl(\frac{5}{6}\Bigr)^{k-1}\!\!=
\boldsymbol{6}$ und $\displaystyle \mathrm{Var}(X_1)=\!\frac{1}{6}\!
\sum_{k=1}^{\infty}k^2\Bigl(\frac{5}{6}\Bigr)^{k-1}\!\!-\!6^2
=\!\boldsymbol{30}$.

Immerhin sind einfache unendliche Reihen (abgeleitete geometrische Reihen) zu
summieren.

Im Allgemeinfall $X_n$ ist ein solches Ausrechnen
"`zu Fu{\ss}"'  m\"uhsam. Elegan\-ter und einfacher ist es,
\textit{bedingte Erwartungswerte} ins Spiel zubringen.

Wir erinnern kurz an dieses Konzept. Ist $F_X$ die
Verteilungsfunktion einer Zu\-falls\-gr\"o{\ss}e und $B$ ein
Ereignis mit $\mathrm{P}(B)>0$, so ist $F_X(x \mid B):= \mathrm{P}(X\leq x \mid B)$ 
die \textit{bedingte Verteilung} und
\[
\boldsymbol{\mathrm{E}(X \mid B):=
\int_{-\infty}^{\infty} x \,\mathrm{d}F_X(x \mid B)}
\]
die \textit{bedingte Erwartung} der Zufallsgr\"o{\ss}e unter der Bedingung
$B$. Im Falle eines diskreten $X$ also einfach 
$\mathrm{E}(X \mid B)= \sum x_k \cdot\mathrm{P}(X=x_k \mid B)$.

Analog zur Formel von der totalen Wahrscheinlichkeit gilt f\"ur eine
disjunkte Zer\-legung $\Omega=B_1\cup B_2\cup \cdots\cup B_n$
mit $\displaystyle \mathrm{P}(B_i)>0$ die
\textit{Erwartungswert-Zerlegungsformel}
\boldmath\[ 
\mathrm{E}(X)=\mathrm{E}(X \mid B_1)
\cdot \mathrm{P}(B_1)+ \cdots+
 \mathrm{E}(X \mid B_n)\cdot \mathrm{P}(B_n).\]\unboldmath%
Diese Formel nutzen wir im Falle $X_n$. Wir betrachten folgende
Ereignisse:

$\boldsymbol{B_1}$: "`Nicht-Sechs im ersten Wurf"';

$\boldsymbol{B_2}$: "`Sechs im ersten und Nicht-Sechs
im zweiten Wurf"';

$\boldsymbol{B_3}$: "`Sechs im ersten und zweiten und
Nicht-Sechs im dritten Wurf"';

$
\phantom{\ \ }\vdots
$

$\boldsymbol{B_n}$: "`Sechs im ersten bis
$(n-1)$-ten und Nicht-Sechs im $n$-ten Wurf"'.

Es gilt \ $\displaystyle \mathrm{P}(B_k)=\frac{5}{6^k}$ \ 
 und \ $\displaystyle \mathrm{E}(X_n \mid B_k)=\mathrm{E}(X_n+k)$ 
\ f\"ur \ $\displaystyle 1\leq k\leq n$. \ \  Ferner

$\boldsymbol{B_{n+1}}$: "`Sechs im ersten bis
$\displaystyle n$-ten Wurf"'; \ \ $\displaystyle
\mathrm{P}(B_{n+1})=\frac{1}{6^n}$, \ \ $\displaystyle
\mathrm{E}(X_n \mid B_{n+1})=n$.

Man pr\"uft leicht nach (endliche geometrische Reihe):
$\sum_{k=1}^{n+1}\mathrm{P}(B_k)=1$.

Mit der Erwartungswert-Zerlegungsformel folgt f\"ur $a:=\mathrm{E}(X_n)$:

$\displaystyle a=\sum_{k=1}^n\frac{5}{6^k}\cdot(a+k)+\frac{n}{6^n}$ \ und
damit \ $\displaystyle \frac{a}{6^n}=5
\sum_{k=1}^n\frac{k}{6^k}+\frac{n}{6^n}$ , also
\boldmath\[
\mathrm{E}(X_n)=\frac{6}{5}\bigl(6^n-1\bigr)=6+6^2+\cdots+6^n.
\]\unboldmath%
Eine Tabelle:
\begin{center}
\begin{tabular*}{\linewidth}{@{\extracolsep\fill}|c|ccccccc|}
\hline
$\boldsymbol{n}$ & {$1$} & {$2$} & {$3$} & {$4$} & {$5$} & {$6$} & {$7$}\\
\hline
{$\boldsymbol{\mathrm{E}(X_n)}$} &
{$6$} & {$42$} & {$258$} & {$1554$} & {$9330$} & {$55986$} & {$335922$}\\
\hline
\end{tabular*}
\end{center}
Nun zur Varianz. Es gilt $\displaystyle \mathrm{E}(X_n^2 \mid
B_k)=\mathrm{E}((X_n\!+\!k)^2)=
\mathrm{E}(X_n^2)\!+\!\frac{12 k}{5}(6^n\!-\!1)\!+\!k^2$ f\"ur $
1\!\leq \!k\!\leq\! n$ und $\mathrm{E}(X_n^2 \mid
B_{n+1})=n^2$. Also mit $b:=\mathrm{E}(X_n^2)$ 

$\displaystyle b=\Bigl(1-\frac{1}{6^n}\Bigr)b+2\bigl(6^n-1\bigr)
\sum_{k=1}^n\frac{k}{6^{k-1}}+
\sum_{k=1}^n\frac{5 k^2}{6^k}+\frac{n^2}{6^n}$, \ und mit \ $\mathrm{Var}(X_n)=b-a^2$ folgt
\boldmath\[
\mathrm{Var}(X_n)=
\frac{6^{n+1}}{25}\bigl(6^{n+1}-10 n-5\bigr)-\frac{6}{25}=
\bigl(\mathrm{E}(X_n)\bigr)^2-\frac{6}{25}\bigl((10 n-7)6^n+7\bigr).
\]\unboldmath%
{\small Bem.: (i) Standardabweichung wie Mittelwert betragen ca. $
1.2\cdot6^n$. \ \ \ (ii) Die Varianz ist wie der Mit\-telwert stets
\textit{ganzzahlig}, da $(10 n-7)6^n+7$ f\"ur alle
$n\in\mathbb{N}$ Vielfaches von $25$ ist (simple
Induktion).\\
(iii) Die \textit{Existenz} der Erwartungswerte wird bei obiger Argumentation
implizit \textit{vorausgesetzt.} F\"ur die Anzahl $A_m$ der
Wurf\-folgen mit $X_n=m$ gilt $A_m=5(A_{m-1}+
A_{m-2}+\cdots+A_{m-n})$; das zeigt eine Sortierung gem\"a{\ss}
$B_1$, $B_2$, $\ldots$, $B_n$. Mit $p_m=A_m/6^m$ also $
p_m=(5/6)p_{m-1}+(5/36)p_{m-2}+\cdots+(5/6^n)p_{m-n}$. Alle Wurzeln des
charakteristischen Polynoms die\-ser Rekursion liegen innerhalb des Kreises
\ $|z|\!<\!1$, da \ $\displaystyle
\frac{5}{6}\!+\cdots+\!\frac{5}{6^n}\!<\!1$. Daraus folgt die Existenz der
Erwartungswerte.}

\bigskip

Wir rechnen $\mathrm{E}(X_2)$ zum Vergleich
auch noch "`zu Fu{\ss}"' aus:

Dazu ist zu \"uberlegen, welche 0-1-Sequenzen der L\"angen 2, 3,
{\dots} es gibt, die am Schluss 11 und ansonsten nur
\textit{isolierte} Einsen enthalten. Einer
\textit{vor} dem Ende 11 vor\-kommenden 1 folgt
\textit{immer} eine 0. D.h.:\\ 
Die Wurf-Folge setzt sich
zusammen aus einer gewissen Anzahl einzelner Nullen, einer gewissen Anzahl
von Paaren 10, gefolgt am Schluss von 11.\\
 F\"ugen wir $
j$ Nullen und $ k$ Paare 10
\textit{in unterschiedlicher Reihen\-folge} aneinander,
entstehen \textit{verschiedene}
Gesamt\-sequenzen, da, vom Beginn der
Wurf-Fol\-ge an gerechnet, die erste Abweichung voneinander eine 0 bei der
einen, eine 1 bei der anderen Sequenz bedeutet. Die L\"ange jeder solchen
Sequenz, 11 am Schluss mit\-gerechnet, ist jeweils $l=j+2k+2$, ihre Eintretens-Wahrscheinlichkeit betr\"agt
$\displaystyle\Bigl(\frac{5}{6}\Bigr)^j\cdot\Bigl(\frac{5}{36}\Bigr)^k\cdot\frac{1}{36}
=\frac{5^{l-k-2}}{6^l}$, \ und es gibt $\displaystyle
\frac{(j\!+\!k)!}{j!\cdot k!}=\binom{l\!-\!k\!-\!2}{k}$ verschiedene
solche Sequenzen. Damit
\boldmath\[
\mathrm{E}(X_2)=
\sum_{l=2}^{\infty} l\cdot\sum_{k=0}^{\lfloor l/2\rfloor-1}
\binom{l\!-\!k\!-\!2}{k}\cdot\frac{5^{l-k-2}}{6^l}.
\]\unboldmath%
Die Auswertung dieser Reihe scheint auf den ersten Blick schwierig,
gelingt aber mittels folgender Iden\-tit\"at (siehe z.B.
Graham/Knuth/Patashnik, \textit{Concrete Mathematics}, 2nd
ed., p.204, Formel (5.74) ):
\[
S_n(z):=\sum_{k=0}^n\binom{n\!-\!k}{k}\,z^k=\frac{1}{\sqrt{1+4z}}
\left(\Bigl(\,\frac{1+\sqrt{1+4z}}{2}\,\Bigr)^{n+1}-
\Bigl(\,\frac{1-\sqrt{1+4z}}{2}\,\Bigr)^{n+1}\right).
\]
Die Herleitung dieser Iden\-tit\"at bei G/K/P ist verwickelt. Man
kann sie aber auch auf ganz einfache Art beweisen:\\ Mittels der Grundformel
des Pascal\-schen Drei\-ecks, $\displaystyle
\binom{n}{k}+\binom{n}{k+1}=\binom{n+1}{k+1}$, zeigt man leicht
$S_{n+2}(z)=z S_n(z)+S_{n+1}(z)$. Die L\"o\-sung
dieser \textit{linearen Differenzen\-gleichung} ergibt die
Iden\-tit\"at. \ Damit
\[
\mathrm{E}(X_2)=\frac{\sqrt{5}}{90}\sum_{l=2}^{\infty}l\left(\Bigl(
\frac{5\!+\!3\sqrt{5}}{12}\Bigr)^{l-1}-\Bigl(
\frac{5\!-\!3\sqrt{5}}{12}\Bigr)^{l-1}\right)=\frac{(1\!-\!q_1)^{-2}-(1\!-\!q_2)^{-2}}{90/\sqrt{5}}=\boldsymbol{42},
\]
wobei $\displaystyle q_{1,2}=\frac{5\pm3\sqrt{5}}{12}$. ($ \mathrm{Var}(X_2)=1650$ und damit
$ \sigma(X_2)\approx40.62$ ergibt sich durch einfache
Rechnung ohne Zusatz\"uberlegungen.)

\textit{Fazit:} Der Weg "`zu Fu{\ss}"' zu
$\mathrm{E}(X_2)$ ist \textit{viel}
steiniger.

\textit{Aber} diese Methode  liefert \textit{mehr
Informationen} \"uber die Zufalls\-gr\"o{\ss}e, n\"am\-lich die genaue
Verteilung der Wahrscheinlichkeiten auf die ein\-zel\-nen Wurf-Anzahlen und
damit z.B. auch den \textit{Median}:
\[\begin{split}
 \frac{\sqrt{5}}{90}\sum_{l=2}^{29}\left(\Bigl(
\frac{5\!+\!3\sqrt{5}}{12}\Bigr)^{l-1}-\Bigl(
\frac{5\!-\!3\sqrt{5}}{12}\Bigr)^{l-1}\right)&\approx0.49961, 
\\   \frac{\sqrt{5}}{90}\sum_{l=2}^{30}\left(\Bigl(
\frac{5\!+\!3\sqrt{5}}{12}\Bigr)^{l-1}-\Bigl(
\frac{5\!-\!3\sqrt{5}}{12}\Bigr)^{l-1}\right)&\approx0.51178
\end{split}\]
In 50\% der F\"alle braucht man also h\"ochstens $\boldsymbol{29}$ W\"urfe 
bis zur ersten Doppel\-sechs, wenn auch der
Durchschnittswert bei $ 42$ liegt. Analog reichen beim
W\"urfeln bis zur ersten Sechs in 51.8\% der F\"alle $
\boldsymbol{4}$ W\"urfe, obwohl der Durch\-schnitt $6$ betr\"agt.

{\small Hinweis: Den Erwartungswert der Anzahl der Versuche bis zum Eintritt eines
Doppelerfolges bestimmt auch N. Henze, \textit{Stochastik f\"ur Einsteiger,
10. Aufl.}, p. 208, mit der bedingten Erwar\-tung. Dies gab mir die
Anregung, das hier darge\-stell\-te Beispiel II. b) auszuarbeiten.}

\clearpage

{\centering\bfseries\large
III. Kartenverteilungen beim Doppelkopf\par}

Da beim Doppelkopf-Spiel jede Karte \textit{doppelt} auftritt, ist zwischen
den Anzahlen ver\-schiedener tats\"achlicher Kartenverteilungen und den
Anzahlen \textit{logischer} Kar\-tenverteilungen (bei denen es nur auf die
Kartenkonstellationen aus der je\-wei\-ligen Spieler-Perspektive ankommt)
sorgf\"altig zu unterscheiden. 

F\"ur Wahrscheinlichkeitsberechnungen muss man die tats\"achlichen einzelnen
Kar\-ten z\"ahlen, was einfach ist. \ Die Erfassung der \textit{logischen}
Kartenver\-tei\-lun\-gen, um die es hier vor allem gehen soll, erweist sich
als erheblich aufwen\-di\-ger. \ Resultat: W\"ahrend die Anzahl
\textit{aller} Verteilungen $\displaystyle 2.36 \cdot 10^{26}$ betr\"agt,
gibt es nur $\displaystyle 2.25\cdot 10^{21}$ \textit{logisch
verschie\-de\-ne} Spielanfangs-Kartenverteilungen.

\textbf{III. a) \ } Wie wahrscheinlich ist eine
"`Hochzeit"'?

Es gibt insgesamt 48 Karten, von denen jeder der vier Spieler 12 er\-h\"alt.
Es sind die Karten \textit{9, 10, Bube, Dame, K\"onig, As} in den vier
Far\-ben \textit{Kreuz, Pik, Herz} und \textit{Karo} aus zwei klassischen
Skat-Spielen; d.h.: Jede Karte kommt genau \textit{zweimal} vor.

Die h\"ochsten Tr\"umpfe -- abgesehen von den beiden "`Dullen"',
den Herz-Zehnen -- sind die Damen, deren wiederum h\"ochste, die bei\-den
Kreuz-Damen (die "`Re-Damen"'), entscheiden, welche
Zweier-Teams gegeneinander spielen: die Spie\-ler mit den beiden Re-Damen
spielen zusammen gegen die anderen beiden. Er\-h\"alt ein Spieler
\textit{bei\-de} Re-Damen (eine "`Hochzeit"'), kann er
entscheiden, ob er heim\-lich allein gegen die anderen drei spielt
(\textit{"`stille"'} Hochzeit) oder aber einen Mit\-streiter
w\"ahlt; z.B.: "`Der erste Stich in fremder Hand geht mit."'
(Zu Beginn des Spiels ist meistens nicht erkenn\-bar, welche zwei Spieler
jeweils ein Team bil\-den, und dies kann manch\-mal tak\-tisch geschickt
genutzt werden.)

Wir bestimmen nun erst einmal \textit{die Anzahl aller verschiedenen
Kar\-ten\-ver\-tei\-lun\-gen} und dann die Anzahl derjenigen unter ih\-nen,
bei denen \textit{einer} der Spie\-ler \textit{bei\-de} Re-Damen erh\"alt. 

Beim Vorgang des Mischens und Karten-Austeilens ist jede \textit{reale Karte}
ein Indi\-vi\-duum. Es spielt \ keine \ Rolle, \ dass \ jeder
\ Karten\textit{wert} (wie Kreuz Bube oder Herz Zehn) doppelt vorhanden ist.
Es gibt also
\boldmath\[
N=\binom{48}{12}\cdot\binom{36}{12}\cdot\binom{24}{12}
\cdot\binom{12}{12}=235809301462142612780721600 
\approx 2.358\cdot 10^{26}
\]\unboldmath%
verschiedene Kartenverteilungen.

Bei wie vielen Kartenverteilungen erh\"alt \textit{Spieler 1} beide Re-Damen?

Diese Anzahl ist offenbar $\displaystyle
\,\binom{2}{2}\cdot\binom{46}{10}\cdot\binom{36}{12}\cdot\binom{24}{12}
\cdot\binom{12}{12}$. Als Eintretens\-wahr\-schein\-lich\-keit ergibt sich
also $\displaystyle w=\binom{46}{10}/\binom{48}{12}=\frac{11}{188}$. 

Die Wahrscheinlichkeit, dass \textit{irgendein} Spieler bei\-de Re-Da\-men
erh\"alt, ist -- da bei allen vier Spielern dieselbe -- gegeben durch
\boldmath\[
w_{\mathrm{H}}=4\cdot\frac{11}{188}=\frac{11}{47}\approx 0.234
.\]\unboldmath%
\fbox{\parbox{0.975\linewidth}{Bei \textit{ann\"ahernd jedem vierten Spiel} landen beide Kreuz-Damen bei ein und demselben Spieler, d.h es tritt eine \textit{Hochzeit}
auf.}}

Nun die kombinatorisch kniffligere Frage:

\textbf{III. b) \ } Wie viele \textit{logisch verschiedene} Kartenverteilungen
treten beim Dop\-pel\-kopf auf?  Anders gesagt: Wie viele verschiedene
"`Spiele"' gibt es?

Die Antwort auf diese Frage, die mir einst ein Student stellte, ohne dass ich
sie auf Anhieb beantworten konnte, verlangt einige
Fallunter\-schei\-dun\-gen. 

Wir klassifizieren die Kartenverteilungen nach der Anzahl $
\boldsymbol{p}$ der Kar\-ten\-\textit{paare} (zwei gleiche Karten in der
Hand \textit{eines} Spielers sind ein solches Paar) und der Anzahl
$\boldsymbol{s}$ der \textit{Solo}karten (Karten, bei denen
die beiden lo\-gisch gleichen Exem\-plare bei \textit{ver\-schie\-de\-nen}
Spielern gelandet sind). Es gilt stets \ $2 p + s=48$. 

Im Falle $p=24$, d.h. $s=0$, erh\"alt jeder der
Spieler 6 Paa\-re; es gibt also $\displaystyle
N_{24}=\binom{24}{6}\cdot\binom{18}{6}\cdot\binom{12}{6}\cdot
\binom{6}{6}=2308743493056\approx2.3\cdot10^{12}$ sol\-che
Ver\-tei\-lun\-gen.

Im anderen Extremfall $p=0$ und damit $s=48$ erh\"alt jeder
Spieler nur Solo\-kar\-ten. 

F\"ur den ersten Spieler gibt es $\displaystyle \binom{24}{12}$ verschiedene
M\"oglichkeiten. Auch beim zwei\-ten kann noch aus allen $24$
verschiedenen Karten aus\-gew\"ahlt werden. Dabei gibt es jeweils
$\displaystyle \binom{12}{d}\cdot\binom{12}{12\!-\!d}=\binom{12}{d}^2$
\textit{lo\-gisch ver\-schiedene} M\"oglichkeiten, $ d\in\{0,
1, \ldots ,12\}$ Karten \textit{zum zweiten Mal} und damit $
12-d$ der bei Spieler 1 nicht vorkommenden $12$ Kar\-ten zu
vergeben. Danach sind $24-d$ verschiedene Karten ver\-ge\-ben
und daher $d$ Karten noch \textit{doppelt} vorhanden. Der
dritte Spie\-ler erh\"alt je eine der noch doppelten Karten und dann
irgend\-wel\-che $12-d$ der noch $24-2 d$
einfach vorhandenen; das sind ins\-ge\-samt $\displaystyle
\binom{24}{12}\binom{12}{d}^2\binom{24- 2 d}{12-d}$lo\-gisch verschiedene
Kartenverteilun\-gen. Fol\-gen\-de Grafik verdeutliche diese
\textit{lo\-gi\-sche} Aufteilung auf die ver\-schie\-denen Kartenarten und
Spieler (\textit{nicht} die Ver\-teilung der \textit{tat\-s\"ach\-lichen
Karten!)}, dabei in Klammern jeweils der/die Spie\-ler:\vspace*{-0.2cm}

 \includegraphics[width=17cm,height=4.6cm,clip=true]{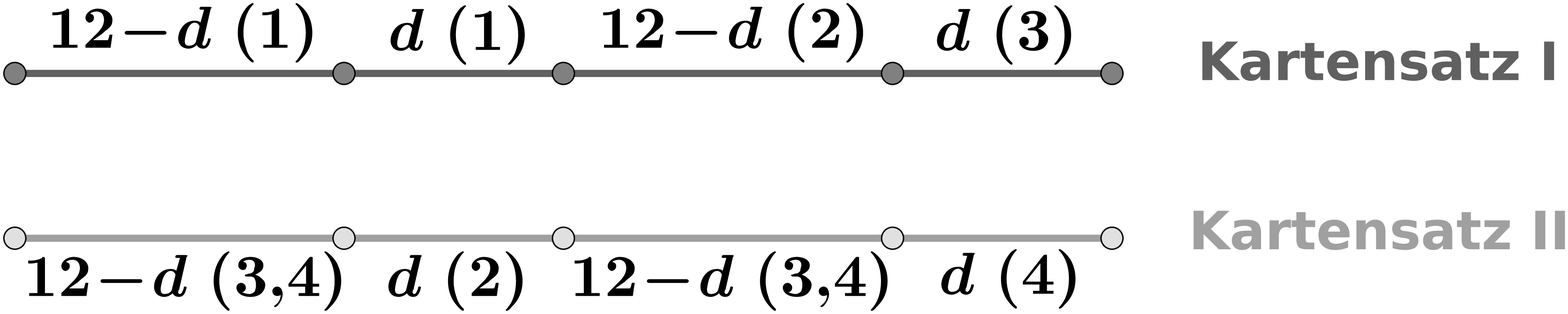}\vspace*{-0.2cm}
 
Die Summe \"uber alle Werte $d$ mit $0 \leq d
\leq 12$ ergibt insgesamt:
\boldmath\[
N_0=\binom{24}{12} \sum_{d=0}^{12} \binom{12}{d}^2\binom{24\!-\!2
d}{12\!-\!d}
= 25780447171287900 \approx 2.578\cdot 10^{16}
\]\unboldmath%
Nun sind noch die Anzahlen $N_p$ der logisch verschiedenen
Karten\-verteilungen f\"ur $0<p<24$ und jeweils $s=48-2 p$ zu bestimmen.

Vorher eine \"Uberlegung zur \textit{Wahrscheinlichkeit }der reinen
Paar-Ver\-teilungen und der reinen Solokarten-Verteilungen:

Bei den Paar-Verteilungen unterscheiden sich logische und tats\"achli\-che
Karten\-ver\-teilungs-Anzahl \textit{nicht}, weil ja beide gleichen Karten
in derselben Hand sind.

Die einzelnen m\"oglichen Paar-Vertei\-lungen sind aus Symmetriegr\"unden
alle un\-ter\-einander gleichwahrschein\-lich. Die Wahrscheinlickeit einer
reinen Paar-Ver\-tei\-lung ist so\-mit
\boldmath\[
w_{\mathrm{P}}=\frac{N_{24}}{N} \approx 9.79\cdot
10^{-15}.
\]\unboldmath%
Die Wahrscheinlichkeit \ f\"ur \textit{einzelne} \ reine Paar-Verteilungen
ist $\boldsymbol{N^{-1} \approx 4.24 \cdot 10^{-27}}$,
entspricht also derjenigen \textit{jeder be\-lie\-bi\-gen} tat\-s\"achlichen
Karten\-ver\-teilung.

Bei den reinen Solokarten-Verteilungen sieht es anders aus. Lo\-gisch
die\-selbe, aber eine faktisch andere Verteilung entsteht, wenn man
ir\-gendwelche der Paa\-re gleicher Karten vertauscht, weil diese dann die
Hand wechseln. Indem man bei gegebener logischer Solo-Ver\-tei\-lung
irgendwie \textit{einen vollen logischen Kar\-ten\-satz} ausw\"ahlt und
die\-sen Karten den "`ersten"' der \textit{faktischen}
Kartens\"atze zuordnet, erh\"alt man insgesamt $\displaystyle 2^{24}$
\textit{verschiedene faktische} Verteilungen zu \textit{ei\-ner
lo\-gi\-schen} So\-lo\-karten-Verteilung. Also sind die
Solokarten-Verteilungen \textit{er\-heb\-lich viel wahrscheinlicher} als die
reinen Paar-Verteilungen. Dies macht schon klar, dass Abz\"ahlen der
\textit{logischen} Verteilun\-gen \textit{nicht} zur Be\-rech\-nung der
\textit{Wahr\-scheinlichkeiten} geeignet ist.

Die Wahrscheinlichkeit einer Solokarten-Verteilung ist
\boldmath\[
w_{\mathrm{S}}=\frac{N_0\cdot 2^{24}}{N} \approx 1.8342\cdot 10^{-3}.\]\unboldmath%
Die Wahrscheinlichkeit jeder \textit{einzelnen} logischen Solokarten-Verteilung ist
dem\-nach

{\centering $\displaystyle
2^{24}/N \approx 7.1147 \cdot 10^{-20}
$.\par}

\fbox{\parbox{0.975\linewidth}{
Die Wahrscheinlichkeit des Auftretens einer reinen Solokarten-Ver\-teilung
be\-tr\"agt immerhin \textit{fast 2 Promille}; eine reine Paar-Ver\-teilung
kommt nur knapp ein\-mal in hundert Billionen Spielen vor, ist also
\textit{\"au{\ss}erst unwahrscheinlich}.}}

Nun m\"ussen zur Bestimmung beliebiger $N_p$ die bisherigen
\"Uberle\-gun\-gen kom\-bi\-niert werden.

Es sind $p$ Kartenpaare und $48-2p$ Solokarten auf die vier Spieler zu verteilen. Es gibt 
$\displaystyle\binom{24}{p}$unterschiedliche M\"oglichkeiten, $p$
Kar\-ten\-paa\-re auszuw\"ahlen. Bei Zu\-ordnung von \ $p_1,
\  p_2$, \ $p_3,  \ p_4$ \ Paaren zu den Spie\-lern \ 1, 2, 3 bzw. 4 \ \ 
(mit $p_1+p_2+p_3+p_4=p$) \ gibt es 
\[
\binom{p}{p_1}\binom{p\!-\!p_1}{p_2}\binom{p\!-\!p_1\!-\!p_2}{p_3}
\binom{p\!-\!p_1\!-\!p_2\!-\!p_3}{p_4}=\frac{p!}{p_1! \,p_2!\, p_3!\, p_4!}=:
\binom{p}{p_1, p_2, p_3, p_4}
\]
lo\-gisch verschiedene Verteilungsm\"oglichkeiten
("`Multinomial-Koeffi\-zi\-ent"').

Au{\ss}erdem erhalten die Spieler noch $12- 2 p_1, 12 - 2 p_2,
12 - 2 p_3 $ bzw. $12-2 p_4$ Solo\-kar\-ten. Daf\"ur stehen
$24 - p$ Karten-Paare zur Ver\-f\"u\-gung. Die Auswahl ist nur
m\"og\-lich, wenn die $p_i$ gewisse Rand\-be\-dingungen
erf\"ullen:

Offenbar ist $\displaystyle p_i \geq \frac{p}{2}-6  \  \Leftrightarrow  \ 
24-p \geq 12 - 2 p_i$ not\-wen\-dig; z.B. im Falle $p=20$ muss
$p_i \geq 4$ gelten. Also
\boldmath\[
\frac{p}{2}-6\leq p_i\leq 6\quad(i=1, 2, 3, 4).\]\unboldmath%

Die \textit{logische} Kartenverteilung verdeutliche wieder ein Diagramm, bei
dem $d_{1\,2}$ die Anzahl der bei den Spielern 1 und 2
\textit{gemeinsam} und $s_1,\  s_2$ die nur bei
\textit{ei\-nem} der Spieler vorkommenden Solokarten be\-zeichne;
$d_{3\,4}, s_3, s_4$ seien analog definiert. 
\begin{center}
\includegraphics[width=17cm,height=4.4cm]{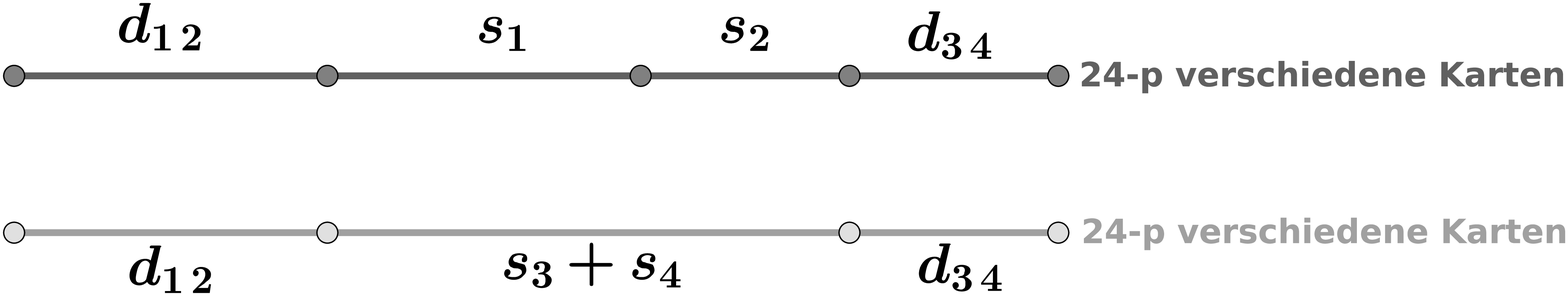}
\end{center}
Es gilt hierbei:
\[
d_{1\,2}+s_1=12-2 p_1, \ d_{1\,2}+s_2=12-2 p_2, \ 
d_{3\,4}+s_3=12-2 p_3, \ d_{3\,4}+s_4=12-2 p_4
\]
sowie
\[
 d_{1\,2}+d_{3\,4}+s_1+s_2=24-p,\quad
d_{1\,2}+d_{3\,4}+s_3+s_4=24-p.
\]
Das sind 6 Gleichungen f\"ur 6 Unbekannte, die allerdings \textit{nicht alle
unabh\"angig} sind: die Summe der ersten vier Gleichungen ist iden\-tisch
mit der Summe der letzten beiden. Legt man aber \textit{einen} der Werte
$d_{1\,2}, d_{3\,4}, s_1, s_2, s_3, s_4 $ fest, ergeben sich
die anderen f\"unf Werte zwangsl\"aufig; das Gleichungssystem hat also den
Rang 5. Legt man z.B. $d_{1\,2}$ auf einen Wert $
\boldsymbol{d}\in\{0, 1, \ldots , \min(12\!-\!2p_1, 12\!-\!2p_2)\}$ fest,
folgt \ $d_{3\,4}=p_1\!+\!p_2\!-\!p_3\!-\!p_4\!+\!d$, und da
dieser Wert nat\"urlich die analoge Bedin\-gung \ $
d_{3\,4}\in\{0, 1, \ldots , \min(12\!-\!2p_3, 12\!-\!2p_4)\}$ erf\"ullen
muss, ergibt sich f\"ur den freien Parameter $d$ insgesamt die
Einschr\"ankung

{\centering
$\displaystyle \max(0,p_3\!+\!p_4\!-\!p_1\!-\!p_2)\leq d \leq \min(12\!-\!2
p_1,12\!-\!2 p_2, 12\!-\!p
\!+\!2 p_3, 12\!-\!p\!+\!2 p_4)$.
\par}

Genau dann, wenn diese Bedingung erf\"ullt ist, sind alle sich zwin\-gend
er\-ge\-ben\-den anderen Werte sinnvoll: 

{\centering
\ $\displaystyle \boldsymbol{
d_{1\,2}=d
}$ \ sowie \ $\displaystyle \boldsymbol{
d_{3\,4} = d\! + \!p_1\! + \!p_2\! - \!p_3 \!- \!p_4, \quad s_1 = 12\! - \!2
p_1 \!- \!d
}$,
\par}

{\centering
$\displaystyle \boldsymbol{
s_2 = 12\! - \!2 p_2\! -\! d, \quad s_3 = 12\! - \!p \!+\! 2 p_4\! - \!d, 
\quad s_4 = 12 \!-\! p\! + \!2 p_3 \!- \!d
}$.
\par}

Nun wird f\"ur jeden zul\"assigen Wert $d$ bei vorgegebenen
$p_1,  \, p_2,  \, p_3 $ und $p_4$ die An\-zahl
der logisch verschiedenen Verteilungen der $48 - 2 p$
Solo\-karten be\-stim\-mt.

Die sich ergebende Anzahl von Solokarten-Auswahlm\"og\-lich\-kei\-ten:
\[
\sum_{d=d_0(p_1, p_2, p_3, p_4)}^{d_1(p_1, p_2, p_3, p_4)}\!\!
\binom{24\! -\!p}{d} \binom{24\!-\!p\!-\!d}{p\!-\!2(p_3\!+\!p_4)\!+\!d}
\binom{24\!-\! 2p_1\!-\! 2 p_2 \!-\! 2 d}{12 \!-\! 2 p_1\! -\! d}
\binom{24\!- \!2 p_1\! -\! 2 p_2\! -\! 2 d}{12\! -\! p \!+\! 2 p_4\! - \!d}
\]
mit
\[
\boldsymbol{d_0(p_1, p_2, p_3, p_4)} = \max(0, p_3+p_4- p_1- p_2)
\]
und
\[
 \boldsymbol{d_1(p_1, p_2, p_3, p_4)} = \min(12 -2p_1, 12-2p_2,
12-p+2p_3, 
12-p+2p_4).
\]
Diese Summe ist dann noch mit der Anzahl der Paar-Auswahlen, also
mit$\displaystyle \binom{24}{p}\binom{p}{p_1, p_2, p_3, p_4}$ zu
multiplizieren.

Wir vergleichen die beiden Grenzf\"alle $p_1=p_2=p_3=p_4=6$ und
$p=0$ mit den vorherigen Formeln (eine Art \textit{Probe}):

Im ersten Fall ergibt sich $0 \leq d \leq 0$ und damit 
\[
\binom{0}{0}\cdot\binom{0}{0}\cdot\binom{0}{0}\cdot\binom{0}{0}\cdot
\binom{24}{24}\cdot\binom{24}{6, 6, 6, 6} =
2308743493056 = N_{24};
\]
im zweiten Fall folgt $0 \leq d \leq 12$ und 
\[
\sum_{d=0}^{12} \binom{24}{d}\binom{24-d}{d}\binom{24 -2 d}{12-d}
\binom{24 -2 d}{12 -d} = 2578447171287900 = N_0.
\]

Die Anzahl \textit{$N_p$ logischer Kartenverteilungen
mit $p$ Paa\-ren} ergibt sich durch Sum\-mieren
\"uber alle zul\"assigen $p_i$-Kombinationen.

Durch die Grenzen $ d_0(p_1, p_2, p_3, p_4)\leq d\leq d_1(p_1,
p_2, p_3, p_4)$ bei der Sum\-mation sind die anfangs genannten
$ p_i$-Einschr\"an\-kun\-gen impli\-zit ber\"ucksichtigt, da
ihr Nicht\-erf\"ulltsein $ d_0>d_1$ zur Folge hat. Die
Aus\-wer\-tung der folgenden Summenformel, f\"ur die Handrechnung zu
kompliziert, \"uberlassen wir -- wie schon die bisherigen
Zah\-len\-rech\-nungen -- dem CA-System MAXIMA. Es gilt
\[
\boxed{\begin{split}
\boldsymbol{N_p}=&\binom{24}{p}
\sum_{\substack{p_i\geq 0\\p_1+\cdots+p_4=p}}
\frac{p!}{p_1! p_2! p_3! p_4!} \,\cdot{}\\
&\sum_{d=d_0(p_1, p_2, p_3, p_4)}^{d_1(p_1, p_2, p_3, p_4)}
\binom{24\! -\!p}{d} \binom{24\!-\!p\!-\!d}{p\!-\!2(p_3\!+\!p_4)\!+\!d}
\binom{24\!-\! 2p_1\!-\! 2 p_2 \!-\! 2 d}{12 \!-\! 2 p_1\! -\! d}
\binom{24\!- \!2 p_1\! -\! 2 p_2\! -\! 2 d}{12\! -\! p \!+\! 2 p_4\! - \!d}
\end{split}}
\]
und die \textit{Anzahl aller logisch verschiedenen Kartenverteilungen} ist

{\centering
\ $\displaystyle \boldsymbol{
N_{\mathrm{L}} = \sum_{p=0}^{24} N_p =
2248575441654260591964 \approx 2.25 \cdot 10^{21}
\approx  N/10^5
}$.
\par}

Da bei einer festen logischen Verteilung mit $\displaystyle p$ Paaren durch
Ver\-tau\-schen der $\displaystyle 24-p$ Solokartenpaare \textit{alle}
faktisch verschiedenen, aber logisch gleichen Ver\-teilungen entstehen, gilt
(eine \textit{Probe!}) \  \ 
$\boldsymbol{\sum_{p=0}^{24}2^{24-p} \cdot N_p=N}\,.
$

\bigskip

\textbf{III. c) \ } Nun beschreiben wir die Zahl der Paare als
\textit{Zufallsgr\"o{\ss}e} $\displaystyle
\boldsymbol{P}$.

Instruktiv ist eine tabellarische \"Ubersicht des Werteverlaufs von
$\displaystyle N_p$ sowie der H\"au\-figkeiten $\displaystyle
H_p=2^{24-p}\,N_p$; \ \ $\displaystyle \boldsymbol{w(P\!=\!p) =
\frac{H_p}{N}}$.

\begin{tabular*}{\linewidth}{@{\extracolsep\fill}|>{$}c<{$}| >{$}r<{$}  >{$}r<{$} |}
\hline
\boldsymbol{p} & \boldsymbol{N_p} & \boldsymbol{H_p=2^{24-p}\cdot N_p}\\
\hline
\hline
0 &25780447171287900=2.58 \cdot 10^{16} & 432524130769286096486400=4.33 \cdot 10^{23}\\ \hline
1 &363832959548298240=3.64 \cdot 10^{17}  & 3052052075130531002449920=3.05\cdot 10^{24} \\ \hline
2 & 2483164491444917280=2.48 \cdot 10^{18}& 10415146759125382327173120=1.04\cdot 10^{25}\\ \hline
3 & 10898289332076706560=1.09 \cdot 10^{19} & 22855369269343329315717120=2.29\cdot 10^{25}\\ \hline
4 & 34502487882513290640=3.45 \cdot 10^{19} & 36178480733894256246128640=3.62\cdot 10^{25}\\ \hline
5 & 83781446916501734400=8.38 \cdot 10^{19}& 43925607240958861325107200=4.39\cdot 10^{25} \\ \hline
6 & 162010440901293880896=1.62 \cdot 10^{20}  &   42470065019628783113601024=4.25\cdot 10^{25} \\ \hline
7 & 255719131625882276352=2.56 \cdot 10^{20}   &  33517618020467641726009344=3.35\cdot 10^{25} \\ \hline
8 & 335095683107207742720=3.35 \cdot 10^{20}   &  21960830688113966626897920=2.20\cdot 10^{25} \\ \hline
9 & 368874473773917855744=3.69 \cdot 10^{20}  &   12087278756623740297019392=1.21\cdot 10^{25} \\ \hline
10 & 343878632431475581440=3.44 \cdot 10^{20}   &  5634107513757295926312960=5.63 \cdot 10^{24} \\ \hline
11 & 272908704913714501632=2.73 \cdot 10^{20}  &   2235668110653149197369344=2.24 \cdot 10^{24} \\ \hline
12 & 184916672029927341504=1.85 \cdot 10^{20}  &   757418688634582390800384=7.57 \cdot 10^{23} \\ \hline
13 & 107050155990705607680=1.07 \cdot 10^{20}  &   219238719468965084528640=2.19 \cdot 10^{23} \\ \hline
14 & 52888925493994157568=5.29 \cdot 10^{19}  &    54158259705850017349632=5.42 \cdot 10^{22} \\ \hline
15 & 22214940975304197120=2.22 \cdot 10^{19}   &   11374049779355748925440=1.14 \cdot 10^{22} \\ \hline
16 & 7897479700512402000=7.90 \cdot 10^{18}    &   2021754803331174912000=2.02 \cdot 10^{21} \\ \hline
17 & 2347160984780451840=2.35 \cdot 10^{18}  &     300436606051897835520=3.00 \cdot 10^{20} \\ \hline
18 & 581872622554903680=5.82 \cdot 10^{17}    &    37239847843513835520=3.72 \cdot 10^{19} \\ \hline
19 & 114698376735022080=1.15 \cdot 10^{17}    &    3670348055520706560=3.67 \cdot 10^{18} \\ \hline
20 & 19220289579691200=1.92 \cdot 10^{16}    &     307524633275059200=3.08 \cdot 10^{17} \\ \hline
21 & 1994754378000384=1.99 \cdot 10^{15}    &      15958035024003072=1.60 \cdot 10^{16} \\ \hline
22 & 249344297250048=2.49 \cdot 10^{14}     &      997377189000192=9.97 \cdot 10^{14} \\ \hline
23 & 0 & 0 \\ \hline
24 & 2308743493056=2.31 \cdot 10^{12}       &      2308743493056=2.31 \cdot 10^{12} \\ \hline
\end{tabular*}

Bemerkenswert: \ $p=23$ ist die einzige An\-zahl an Paaren, die
bei kei\-ner Kar\-ten\-verteilung vorkommt. Was leicht un\-mit\-telbar
einzusehen ist: Ins\-gesamt zwei Solokarten (bei ver\-schie\-denen Spielern
anzutref\-fen) kann es nicht geben, da alle Spieler $12$
Karten erhalten, eine \textit{gera\-de} Anzahl.

94.74\% aller Kartenverteilungen, n"amlich \ $223410396488155960977653760$,
haben $2$
bis $ 9$ Kartenpaare; die \ $156091771014949542410846208$
Verteilungen mit $4$ bis $
7$ Kar\-ten\-paaren sind schon ca. zwei Drittel aller
Verteilun\-gen.

Mehr als $
12$ Kartenpaa\-re haben \ $
287134455041807189648832$, also nur 1.22{\textperthousand} aller
Kartenverteilungen; rei\-ne Solo\-kar\-ten-Ver\-teilungen ($
p=0$) sind etwa eineinhalb mal so h\"aufig.

(Hier ist nicht von den logischen Anzahlen $N_p$, sondern von
den \textit{H\"au\-fig\-keiten} $H_p$ und daher von
\textit{Wahr\-scheinlich\-kei\-ten} die Rede. Die gr\"o{\ss}te Anzahl
\textit{logisch ver\-schie\-dener Spiele} gibt`s bei $9$
Paaren, 90.3\% aller lo\-gisch verschiedenen Spiele haben $6$
bis $13$, 70.1\% haben $7$ bis $
11$ und nur 23.8\% haben $4$ bis $7$ Paare.)

Der \textit{Erwartungswert} der Paar-Anzahl $P$ ist
\ $\displaystyle
\mathrm{E}(P)\!=\!\sum_{p=0}^{24}p\,\frac{H_p}{N}\!=\!\frac{264}{47}
\!\approx\! 5.617$, die \textit{Va\-ri\-anz} be\-tr\"agt $\displaystyle
\mathrm{Var}(P)\!=\!\frac{48576}{11045} \!\approx\! 4.398$ und damit die
\textit{Standardab\-wei\-chung} $\sigma(P)\!\approx\!2.09714$.

Da $\displaystyle \frac{1}{N}\sum_{p=0}^{5} H_p \approx 0.49556645766$ und
$\displaystyle \frac{1}{N}\sum_{p=0}^{6} H_p \approx 0.67566989$, ist
$\displaystyle 5.5$ der \textit{Me\-di\-an} der Paar-Anzahl $\displaystyle
P$. Der \textit{Modus} ist $ 5$, aber $ 6$ ist
fast ge\-nau\-so h\"aufig.

Einen noch an\-schaulicheren Eindruck verschaf\-fen \textit{grafi\-sche}
Darstel\-lun\-gen und der Ver\-gleich mit der
\textit{Binomialver\-tei\-lung} (links) bzw. \textit{Nor\-malver\-teilung}
(rechts, mit Ste\-tig\-keits\-kor\-rek\-tur) mit passenden
Kenn\-gr\"o{\ss}en; es ergeben sich recht gute \"Uber\-ein\-stimmungen:

\includegraphics[width=8.4cm,height=8.3cm,clip=true]{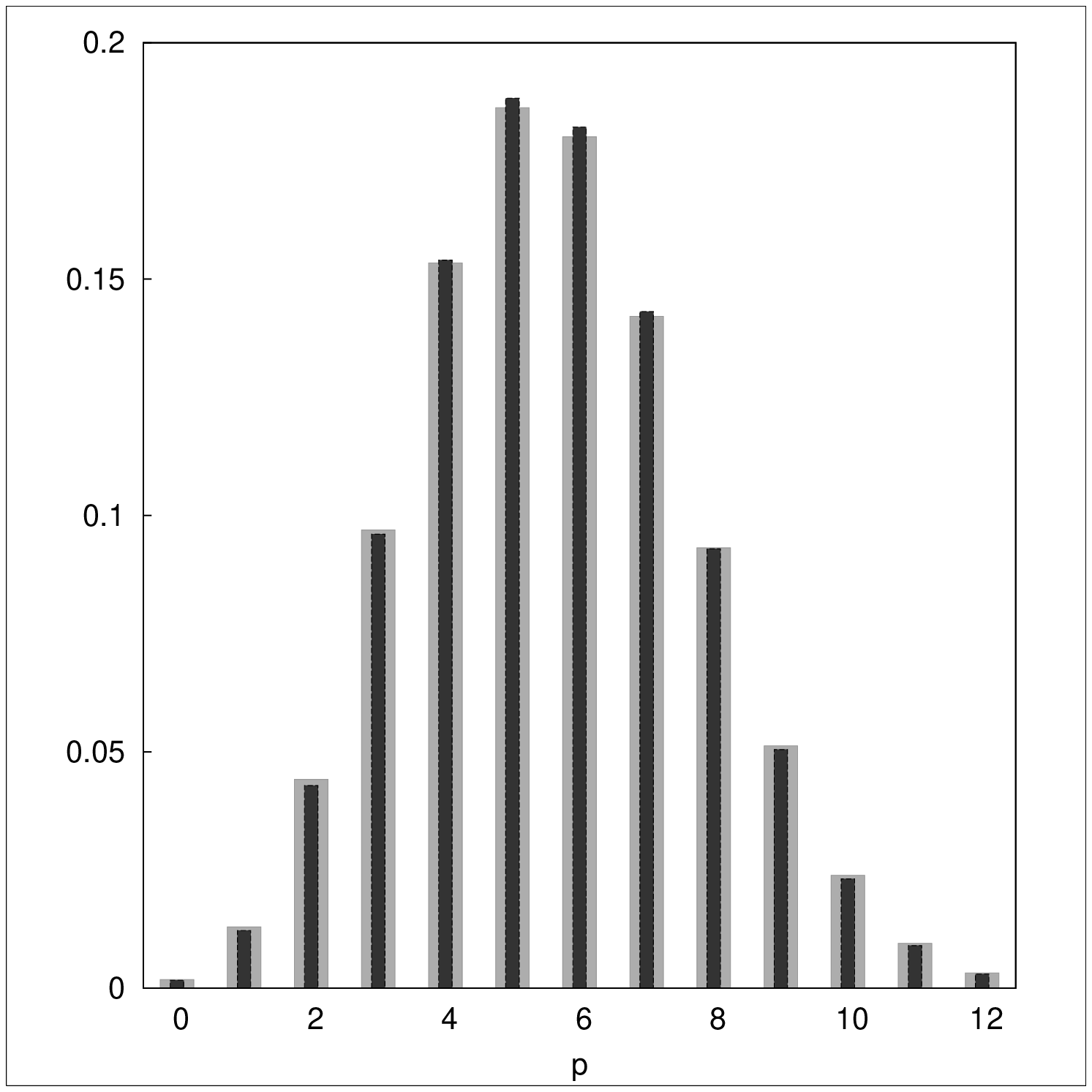}
\     
\includegraphics[width=8.4cm,height=8.3cm,clip=true]{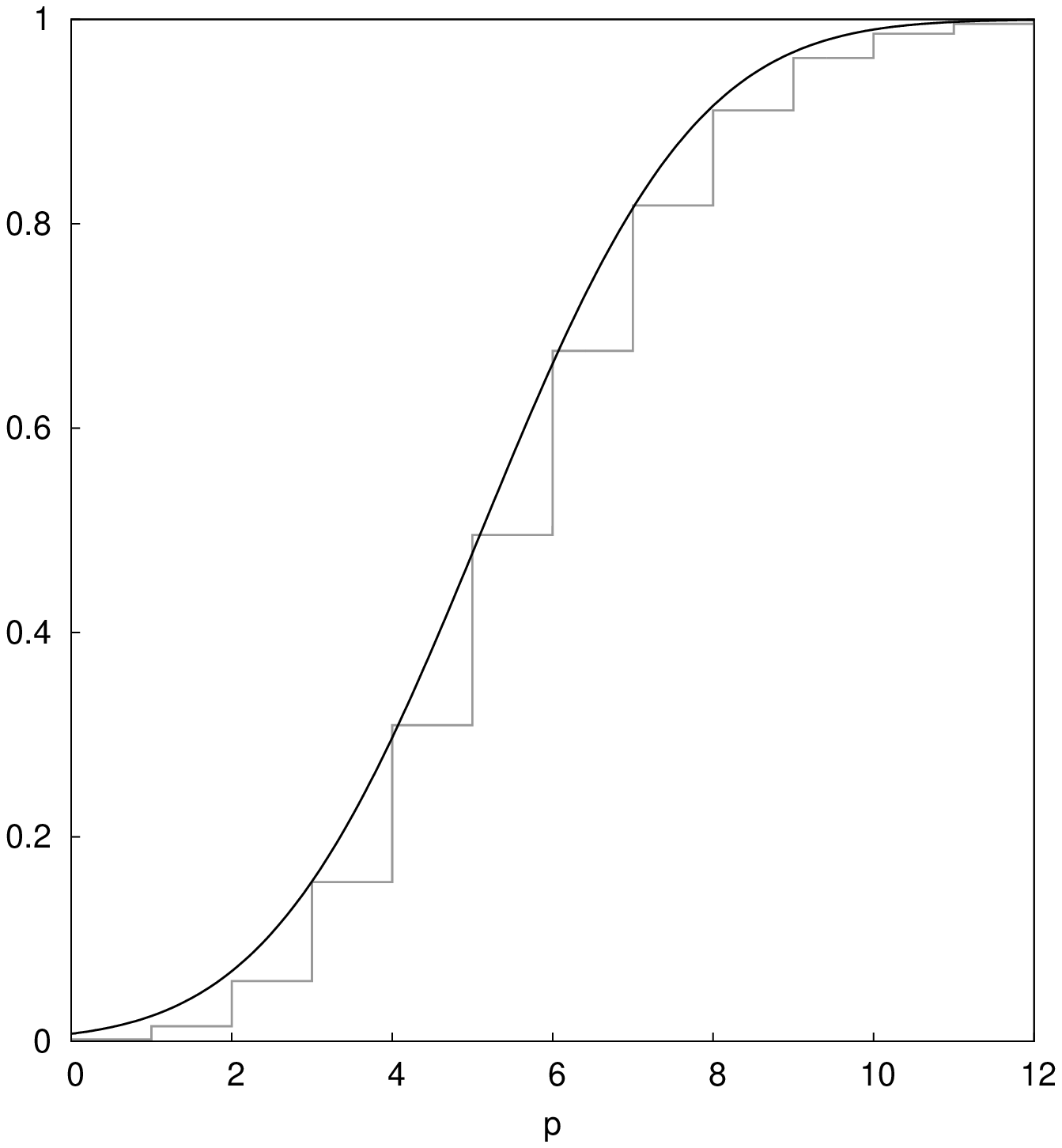}

Die Bi\-no\-mial\-ver\-tei\-lung (schwarz, $\displaystyle n=24,
p=\frac{11}{47}$) hat eine geringf\"ugig klei\-ne\-re Va\-ri\-anz als die
tats\"ach\-li\-che Paar-Verteilung (grau). Die Wahr\-schein\-lich\-keit,
dass ir\-gendein Kartenpaar in \textit{einer} Hand landet, be\-tr\"agt
jeweils $\displaystyle w_{\mathrm{H}}=\frac{11}{47}$; siehe a). W\"aren
diese Er\-eig\-nisse alle \textit{unabh\"angig} voneinan\-der, w\"are die
Binomial\-ver\-tei\-lung die exakte Verteilung der Zufallsgr\"o{\ss}e
$P$; auf jeden Fall \ $\displaystyle
\mathrm{E}(P)=24\cdot\frac{11}{47}$ .

% \bigskip
% 
% {\small Der zentrale Grenzwertsatz, immer formuliert als Aussage \"uber Summen
% \textit{unab\-h\"an\-giger} Zufallsvariablen -- siehe etwa H.~Bauer,
% \textit{Wahrschein\-lich\-keits\-theorie, 5. Auf\-lage}, de Gruy\-ter 2002,
% S. 240f.; A.~R\'enyi, \textit{Probability Theory}, North
% Hol\-land/\-Akad\'emiai Kiad\'o 1970 (Nach\-druck Dover 2007), chap. VIII;
% W.~Feller, \textit{An Intro\-duction to Probability Theory and Its
% Appli\-cations, Volume II}, Wiley 1971, chap. IX --, gilt hier
% n\"aherungsweise auch f\"ur eine endliche Summe nur "`schwach ab\-h\"an\-giger"'
%  Zufallsgr\"o{\ss}en, der 24 Indika\-to\-ren des
% Auftretens der einzelnen Karten\-paare.}

\clearpage

{\centering\bfseries\large
IV. Zufallsprodukte\par}

Seien $X_1, X_2, \ldots , X_n$ unabh\"angige auf dem Intervall
$[0,A]$ mit einem $A>0$ gleich\-verteilte
Zufallszahlen.

Es sollen die Eigenschaften der Zufallsgr\"o{\ss}en $
Y_n:=X_1\cdot X_2\cdots X_n \ (n=1, 2, \ldots)$ un\-ter\-sucht wer\-den.

Vor einigen Jahren experimentierte ein Student der Informatik mit solchen
Zu\-falls\-zahlen auf dem Computer, war \"uberrascht, dass im Falle
$A=2$ die $Y_n$ meist \textit{sehr klein}
ausfielen, obwohl doch ihr Mittelwert $1$ sein sollte, und
sprach mich deshalb an. Erst f\"ur $A=2.72$ (so in etwa)
\"andere sich das seltsamerweise, und f\"ur noch gr\"o{\ss}ere
$A$ fielen die Zufallsprodukte meist \textit{sehr gro{\ss}}
aus.

Ich best\"atigte ihm die mathematische Korrektheit seiner \"Uberlegung, der
Mittelwert sei $(A/2)^n$, also genau im Falle $
A=2$ gleich $1$. Wieso dann dieses seltsame Verhalten der
Zufallsprodukte?

Indem man die Logarithmen $\ln Y_n = \ln X_1 + \cdots +\ln X_n$
betrachtet, erh\"alt man den ers\-ten Hin\-weis darauf, woran es liegt:
$\ln Y_n$ ist ann\"ahernd nor\-mal\-verteilt und hat den
Mit\-telwert \ $n \cdot \int_0^{A} \ln x \,{\small
\frac{1}{A}}\,\mathrm{d}x=n(\ln(A)-1)$. 
Also \ \ \ $0.5\approx\mathrm{P}\bigl(\ln Y_n\leq
n(\ln(A)-1)\bigr)
=\mathrm{P}\bigl(Y_n\leq (A/\mathrm{e})^n\bigr)$. \ \ \ Und das bedeutet: 

Im Falle \ $A=\mathrm{e}\approx2.72$ \ gilt bei etwa 50\% aller
Werte \ $Y_n \leq 1$, und f\"ur die andere H\"alfte gilt
$Y_n\geq 1$; Mittelwert und Medi\-an klaffen f\"ur
gr\"o{\ss}ere $n$ sehr weit aus\-ein\-ander.

(Am Beispiel $\varphi=\ln$ zeigt sich hier, dass im allgemeinen
$
\mathrm{E}\bigl(\varphi(X)\bigr)\neq\varphi\bigl(\mathrm{E}(X)\bigr)$ gilt.)

Dem Studenten reichte diese Ad-hoc-Erkl\"arung. Aber ich wollte es doch etwas
ge\-nauer wissen und habe mich daher ein bisschen mit den Eigenschaften der
Gr\"o{\ss}en $Y_n$ besch\"aftigt, was sich als nicht v\"ollig
trivial erwies.

\bigskip

{\bfseries
IV. a) \ }  Vorbemerkung zu \textit{Summen}
gleichverteilter Zufallsgr\"o{\ss}en

Die Dichte der Summe $X+Y$ zweier
unabh\"angiger stetiger Gr\"o{\ss}en ergibt sich be\-kannt\-lich als das
Faltungsprodukt 
\[
f_{X+Y}(x)=\int_{-\infty}^{\infty}f_X(x-t)f_Y(t)\,\mathrm{d}t.
\]
Benutzt man dies f\"ur unabh\"angige Gleichverteilungen auf dem
Inter\-vall $[-1,1]$, ergibt sich induktiv f\"ur
$n$ Summanden die Formel
\[
f_n(x)=\begin{cases}
\frac{1}{2^n\,(n-1)!}{\displaystyle\sum_{k=0}^{\lfloor(n+x)/2\rfloor}}
(-1)^k\binom{n}{k}(n\!+\!x\!-\!2k)^{n-1}, & |x|<n,\\
0,& \mathrm{sonst},
\end{cases}
\]
die schon N. I. Lobatschewski, bekannt vor allem durch seine
Mitbegr\"undung der nichteuklidischen Geometrie, im Zusammenhang mit
\"Uberlegungen zur Feh\-lerrechnung angegeben hat. (Siehe A. R\'enyi,
\textit{Probability Theory},  p. 197f.)

Hat $X$ die Dichte $f$, so hat $
\tilde{X}=aX+b$ mit $a>0$ bekanntlich die Dichte 
$\displaystyle \tilde{f}(x)=\frac{1}{a}\, f \bigl( \frac{x-b}{a} \bigr)$.

Umgerechnet auf den Fall der Zufallsgr\"o{\ss}en $\displaystyle X_1 , X_2 ,
\ldots, X_n$ erhalten wir mit $\displaystyle a=A/2$, $\displaystyle b=n\cdot
A/2$
\[
f_{X_1 + \cdots+ X_n}(x)=\begin{cases}
\frac{1}{A\cdot(n-1)!}{\displaystyle\sum_{k=0}^{\lfloor\frac{x}{A}\rfloor}}
(-1)^k\binom{n}{k}\bigl(\frac{x}{A}\!-\!k\bigr)^{n-1}, & 0<x<n A,\\
0,& \mathrm{sonst}.
\end{cases}
\]
Als Verteilungsfunktion folgt
\boldmath\[
F_{X_1 + \cdots+ X_n}(x)=\begin{cases}
0, & x<0,\\
\frac{1}{n!}{\displaystyle\sum_{k=0}^{\lfloor\frac{x}{A}\rfloor}}
(-1)^k\binom{n}{k}\bigl(\frac{x}{A}\!-\!k\bigr)^{n}, & 0 \leq x<n A,\\
1,& x\geq n A.
\end{cases}
\]\unboldmath%

{\bfseries
IV. b) \ } Dichte und Verteilung von $Y_n$

Per Induktion k\"onnen wir die Verteilungsfunktion $F_n$ zu
$Y_n$ bestimmen:
\[
 F_1(x)=\mathrm{P}(X_1\leq x)=\int_0^x
 \frac{1}{A}\,\mathrm{d}x_1=\frac{x}{A} \ \  (0\leq x\leq A),
\]
\[\begin{split}
F_2(x) &=\mathrm{P}(X_1\cdot X_2\leq x)=\int_0^A\int_0^{\min(A,\,x/x_1)}
\frac{\mathrm{d}x_2}{A}\,\,\frac{\mathrm{d}x_1}{A}
=\int_0^A \frac{\min(A,x/x_1)}{A^2} \,\mathrm{d}x_1
\\
&=\int_0^{x/A}\frac{\mathrm{d}x_1}{A}
+\int_{x/A}^A\frac{x/x_1}{A^2}\,\mathrm{d}x_1
=\frac{x}{A^2}\bigl(1+(2\ln A -\ln x)\bigr) \ \ (0 \leq x\leq A^2).
\end{split}\]
Schluss von $\displaystyle n$ auf $\displaystyle n+1$ mittels analoger
Umformung ergibt: 
\boldmath\[
\boxed{F_n(x):=F_{Y_n}(x)=\mathrm{P}
(X_1\,X_2 \cdots X_n\!\leq\! x)=\frac{x}{A^n}
\sum_{k=0}^{n-1}\frac{(n \ln A -\ln x)^k}{k!} \ \ (0\!\leq\! x\!\leq \!A^n)}
\]\unboldmath%
Und zwar 
\[\begin{split}
&\mathrm{P}(X_1\cdots X_n\cdot X_{n+1}\leq x) = 
\int_0^A \mathrm{P}\left(X_1\cdots X_n\leq \min \Bigl(
A^n,\frac{x}{x_{n+1}}\Bigr)\right)\frac{\mathrm{d}x_{n+1}}{A}
\\
&= \int_0^{x/A^n}\mathrm{P}(X_1\cdots X_n\leq A^n)\,
\frac{\mathrm{d}x_{n+1}}{A} +
\int_{x/A^n}^A \mathrm{P}\Bigl(X_1\cdots X_n\leq
\frac{x}{x_{n+1}}\Bigr)\,\frac{\mathrm{d}x_{n+1}}{A}
\\
&=\frac{x}{A^{n+1}} + \int_{x/A^n}^A  \frac{x}{A^n\cdot x_{n+1}}\cdot
\sum_{k=0}^{n-1}\frac{(n \ln A -\ln x + \ln x_{n+1})^k}{k!}\;\frac{
  \mathrm{d}x_{n+1}}{A}
\\
&= \frac{x}{A^{n+1}} + \frac{x}{A^{n+1}}\cdot
\left.\sum_{k=0}^{n-1}
\frac{(n \ln A -\ln x + \ln x_{n+1})^{k+1}}{(k+1)!}\right|_{
x_{n+1}=x/A^n}^A
\\
&= \frac{x}{A^{n+1}} + \frac{x}{A^{n+1}}\cdot 
\sum_{k=1}^n \frac{((n+1) \ln A -\ln x)^k}{k!}.
\end{split}\]
Mit der gemeinsamen Dichte $\displaystyle f(x_1 , \ldots,x_n)= 1/A^n \ \ 
(0\leq x_1 , \ldots, x_n\leq A)$ des Zu\-falls\-vektors $\displaystyle (X_1
, \ldots X_n)$ ist es praktisch eine reine \textit{Volumen}berechnung.

Durch Ableiten erhalten wir die zugeh\"origen Dichten $f_n$:
\boldmath\[
f_n(x):=f_{Y_n}(x)=\frac{(n\,\ln A-\ln x)^{n-1}}{A^n\cdot (n-1)!}
\ \ (0 < x  < A^n)
\]\unboldmath%
Man beachte, dass $F_n$ durch den Faktor $x$
stetig, $f_n$ aber bei $0$
\textit{unbeschr\"ankt} ist.

Die folgenden Plots verdeutlichen, wie \textit{klein} im Falle $
A=2$ die meisten Produkte mit wachsendem $ n$ ausfallen; links
$ F_n$ f\"ur $ n=3$, rechts f\"ur $n=24$.

\includegraphics[width=7.4cm,height=5.5cm]{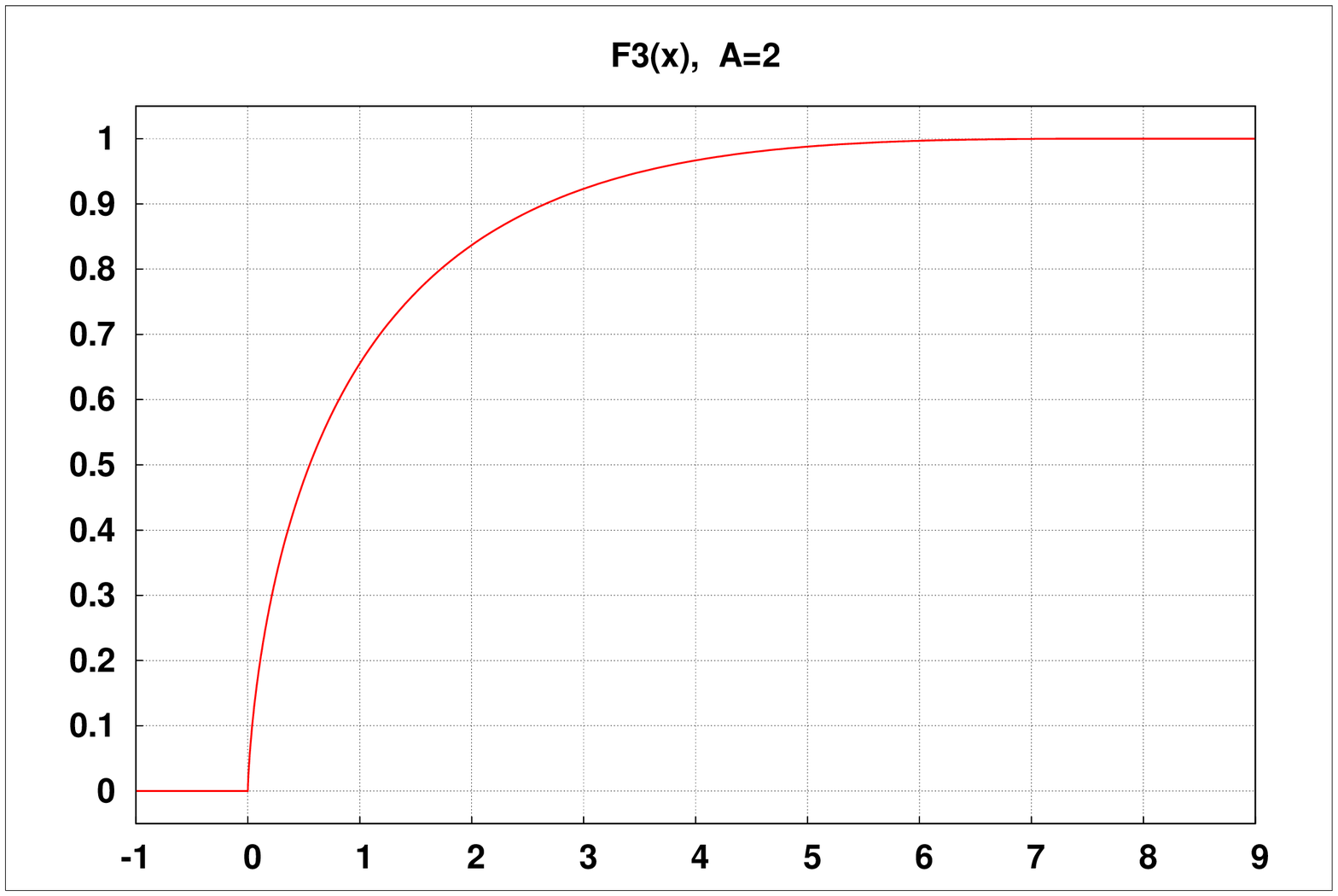}
  \ 
\includegraphics[width=9.4cm,height=5.5cm]{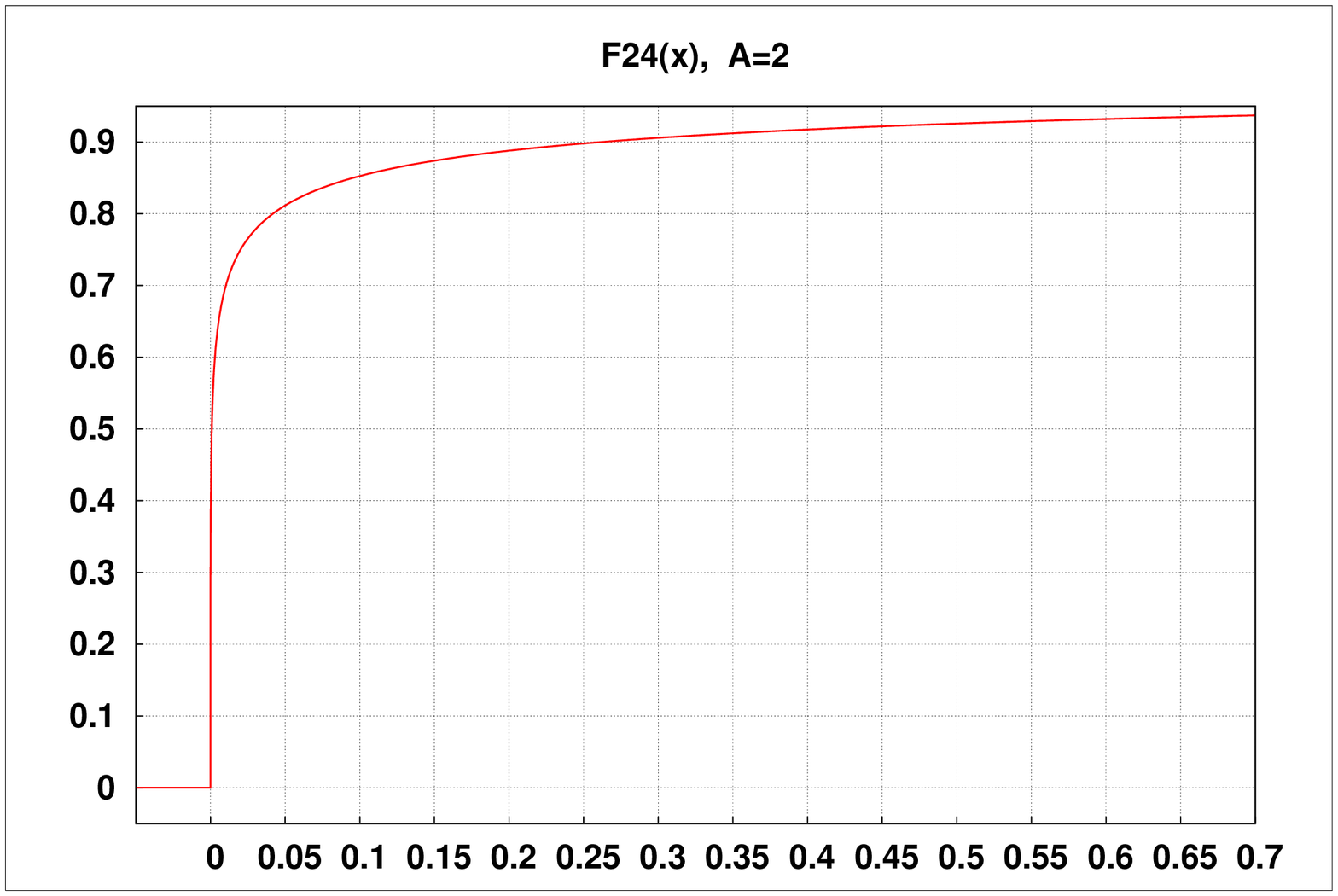}

Die n\"achsten Grafiken zeigen die F\"alle \ $A=\mathrm{e}$
\ und \ $A=1.25\,\mathrm{e}$ ; man sieht den
\textit{gro{\ss}en} Unterschied zwischen beiden schon f\"ur $n=24$.

\includegraphics[width=9.4cm,height=5.5cm]{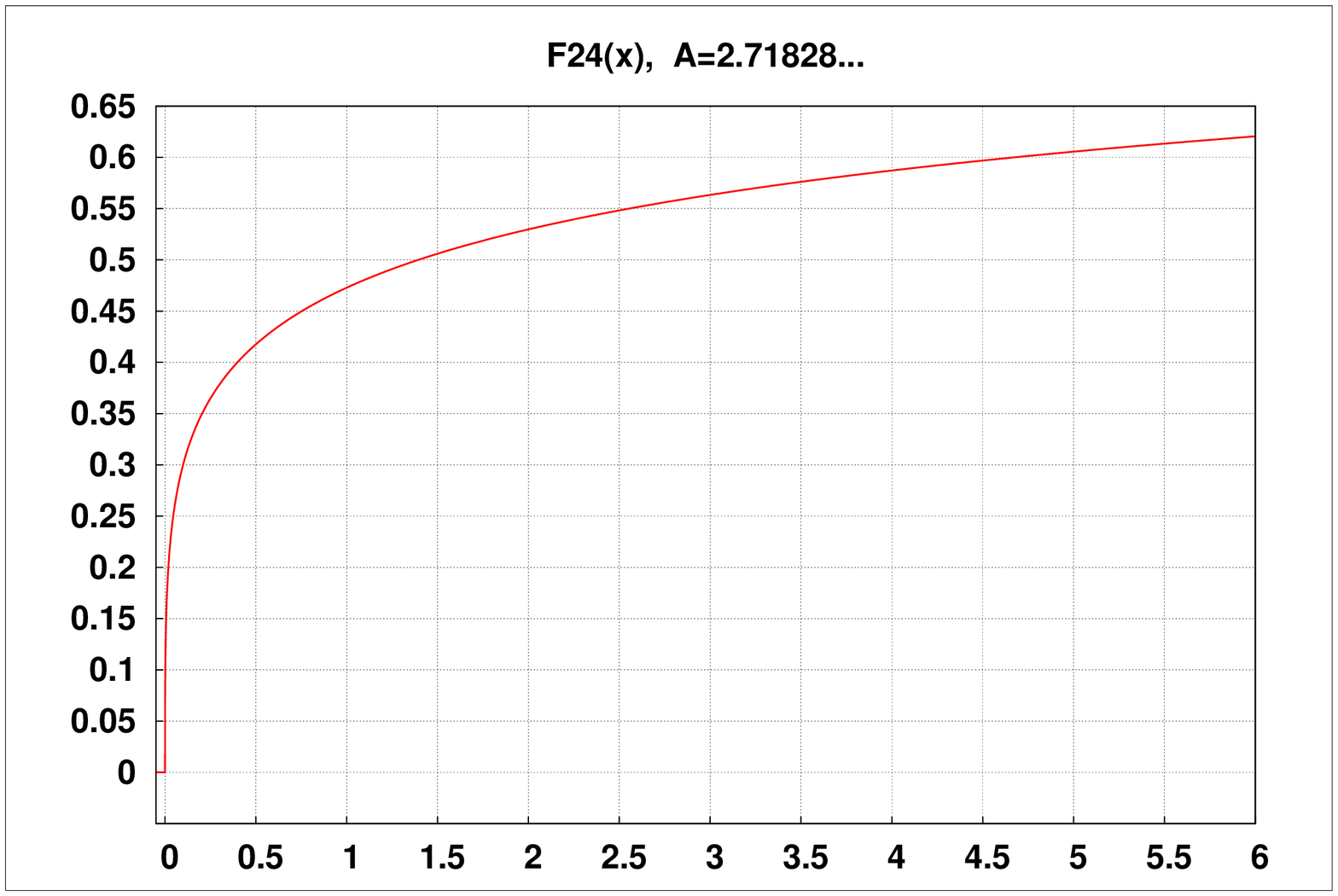}
\ 
\includegraphics[width=7.4cm,height=5.5cm]{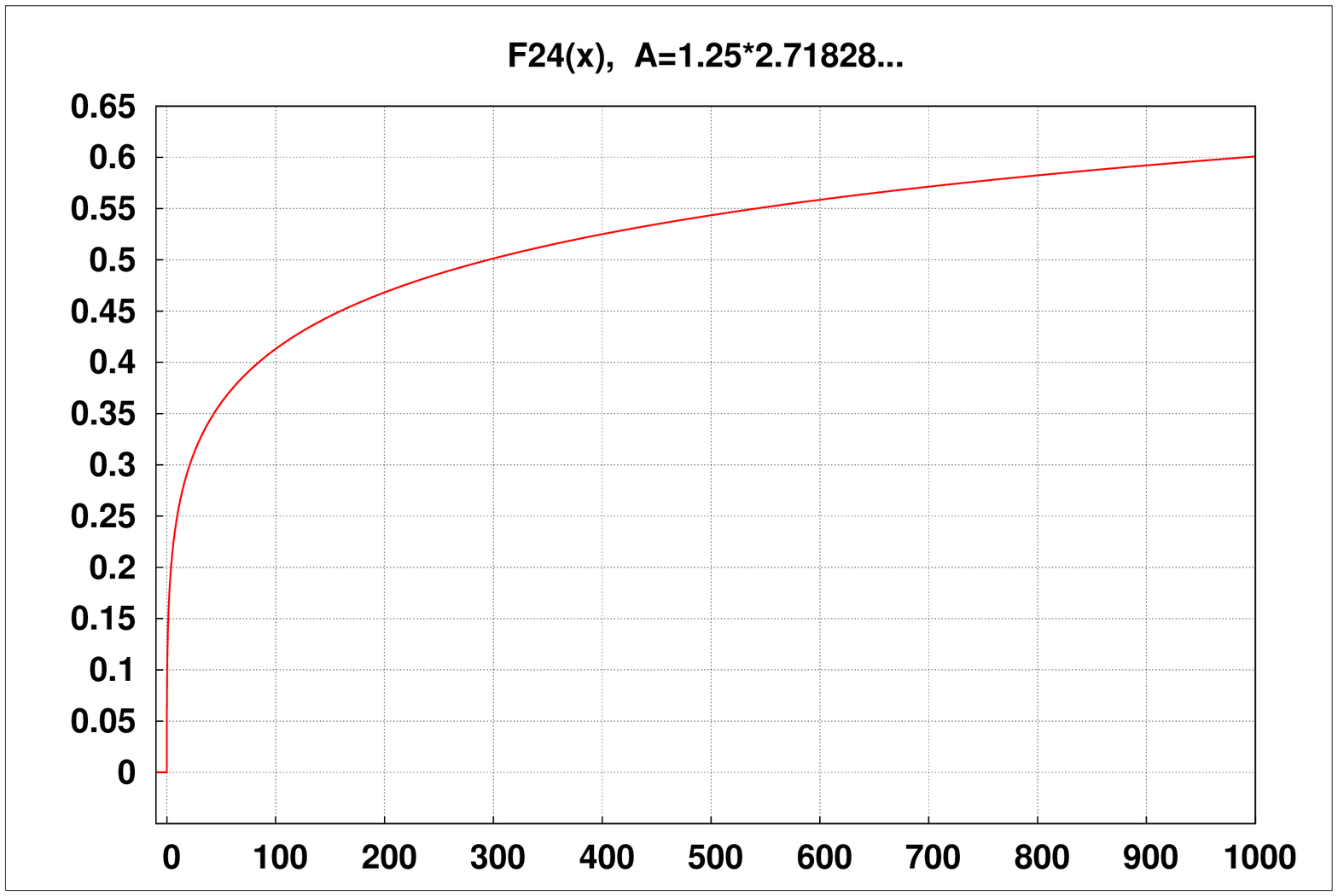}

Die allgemeine Formel f\"ur den \textit{Erwartungswert des Produktes} zweier
Zufalls\-gr\"o{\ss}en lautet
\boldmath
\[
\mathrm{E}(X\,Y)=\mathrm{E}(X)\,\mathrm{E}(Y)+\mathrm{E}\Bigl(
\bigl(X-\mathrm{E}(X)\bigr)\bigl(Y-\mathrm{E}(Y)\bigr)\Bigr);
\]\unboldmath%
im Falle \textit{unabh\"angiger} Gr\"o{\ss}en f\"allt der Kovarianz-Term weg.
Der Beweis beruht auf der Identit\"at $
AB-ab=(A-a)b+a(B-b)+(A-a)(B-b)$ sowie der all\-gemeinen Additivit\"at des
Erwartungswerts. 

Durch analoge elementare, wenn auch etwas umfangreichere Manipulationen
erh\"alt man eine Formel f\"ur die Varianz des Produktes, wobei 
$ \mathrm{Cov}(X,Y):=
\mathrm{E}\Bigl(\bigl(X-\mathrm{E}(X)\bigr)\bigl(Y-\mathrm{E}(Y)\bigr)\Bigr)$
als Notation f\"ur \textit{Kovarianzen} benutzt wird:
\boldmath
\[
\begin{split}\mathrm{Var}(X\,Y) = &\mathrm{Var}(X)\,\mathrm{Var}(Y)+
\mathrm{E}(X)^2\,\mathrm{Var}(Y)+\mathrm{E}(Y)^2\,\mathrm{Var}(X)\\
&+\mathrm{Cov}(X^2,Y^2)-\mathrm{Cov}(X,Y)\bigl(
2\mathrm{E}(X)\,\mathrm{E}(Y)+\mathrm{Cov}(X,Y)\bigr)
\end{split}\]\unboldmath%
Bei \textit{unabh\"angigen} Zufallsgr\"o{\ss}en (genauer gesagt: wenn sowohl
$X$ und $Y$ als auch $X^2$ und
$Y^2$ \textit{unkorreliert} sind) vereinfacht sich die
Varianz-Formel wesent\-lich durch Weg\-fall aller Kovarianzterme.

Die Anwendung dieser Formeln auf die Zufallsgr\"o{\ss}en $Y_n$
ergibt einerseits 
offenbar $
\boldsymbol{\mathrm{E}(Y_n)=\Big(\frac{A}{2}\Bigr)^n}$ und per Induktion \ 
$\displaystyle
\mathrm{Var}(Y_n)=A^{2n}\left(\Bigl(\frac{1}{3}\Bigr)^n\!\!-\!
\Bigl(\frac{1}{4}\Bigr)^n\right)$;
die Standardabweichung ist also 
\boldmath\[
\sigma(Y_n)=\Bigl(\frac{A}{\sqrt{3}}\Bigr)^n
 \,\sqrt{1-\Bigl(\frac{3}{4}\Bigr)^n}.
\]\unboldmath%
{\small Die allgemeine Varianz-Formel habe ich -- seltsamerweise, da ihre Herleitung
nicht schwer ist -- nicht in der Literatur gefunden. \ Vgl. z.B. H. Rinne,
\textit{Taschenbuch der Statistik}, 3. Auflage, Harri Deutsch 2003, S. 222,
wo nur eine "`N\"aherungsformel"'
 angegeben ist.}

\bigskip

\textbf{IV. c) \ } Die Verteilung von $ L_n:=\ln Y_n =\ln X_1 +
\cdots + \ln X_n$

Diese ergibt sich leicht aus der von $\displaystyle Y_n$ und umgekehrt:
\[ G_n(x):=\mathrm{P}(L_n\leq x)=\mathrm{P}(Y_n\leq
\mathrm{e}^x)=\frac{\mathrm{e}^x}{A^n}\sum_{k=0}^{n-1}
\frac{(n \ln A -x)^k}{k!} \ \ (0 < \mathrm{e}^x<A^n),
\]
und die Dichte, also die Ableitung, ist \ $\displaystyle
g_n(x)=\frac{\mathrm{e}^x\,(n \ln A-x)^{n-1}}{A^n\,(n-1)!}
 \ \,(-\infty<x< n \ln A)$.

Plots zeigen, dass die Dichten $g_n$ wirklich an\-n\"ahernd
nor\-mal\-verteilt aussehen, f\"ur kleinere $n$ nat\"urlich
noch unvollkommen; links $n=12$, rechts $n=120$,
\ $A=2$:

\includegraphics[width=8.4cm,height=5.5cm]{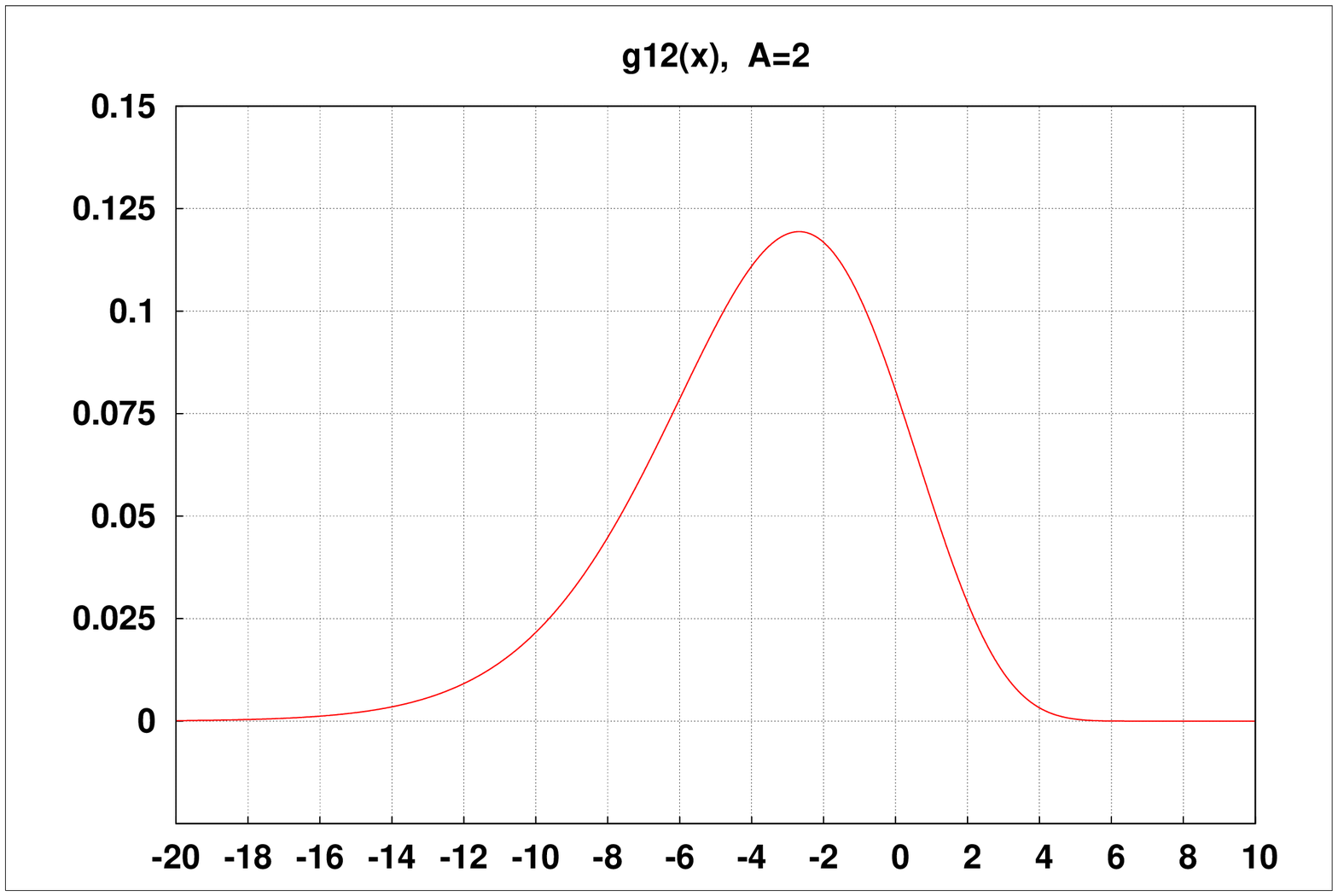}
 \ 
\includegraphics[width=8.4cm,height=5.5cm]{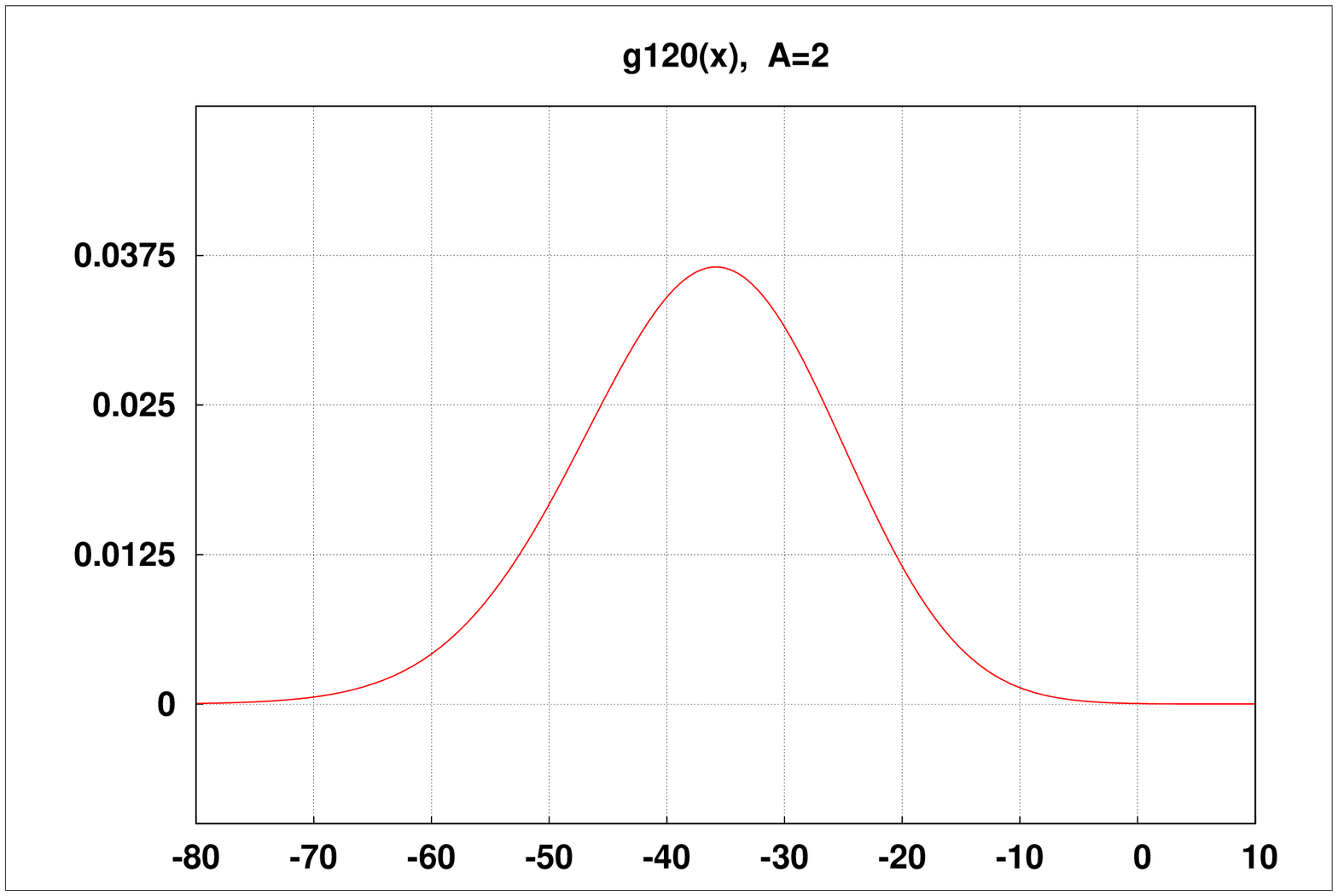}

\bigskip

Man kann die Dichten $g_n$ auch \textit{unmittelbar},
unabh\"angig von $F_n$, \"uber \textit{Fourier\-integrale}
bestimmen:

Ist $f_X$ die Dichte einer stetigen Zufallsgr\"o{\ss}e
$X$, ist 
\[
\varphi_X(t):=\int_{-\infty}^{\infty}\mathrm{e}^{\mathrm{i}tx}
f_X(x)\,\mathrm{d}x=\mathrm{E}\bigl(\mathrm{e}^{\mathrm{i}tX}\bigr)
\]
die \textit{charakteristische Funktion} von $X$. \ \ Im
Allgemeinfall $\displaystyle
\varphi_X(t):=\int_{-\infty}^{\infty}\mathrm{e}^{\mathrm{i}tx}\,
\mathrm{d}F_X(x)$. Per Fourier-Umkehrtheorem sind dann auch $
f_X$ und die Verteilungs\-funk\-tion $F_X^{\mathstrut}$ durch
$\varphi_X$ darstellbar, falls $\varphi_X$
absolut integrierbar:
\[
f_X(x)=\frac{1}{2\pi}\int_{-\infty}^{\infty}\mathrm{e}^{-\mathrm{i}tx}\varphi_X(t)\,\mathrm{d}t.
\]
Allgemeinere Umkehraussagen hinsichtlich der \textit{Verteilungsfunktion}
findet man z.B. bei H. Cram\'er, \textit{Mathematical Methods of
Statistics}, \ Princeton 1946, $^{19}$1999.

Da die Dichte einer unabh\"angigen Summe $X+Y$ durch die
Faltung der Dichten $f_X$ und $f_Y $ gegeben
ist, gilt gem\"a{\ss} dem Fourier-Faltungs\-satz \ die
\textit{Produktformel}
\boldmath\[\varphi_{X+Y}(t)=\varphi_X(t)\cdot\varphi_Y(t).
\]\unboldmath%
Mit \ $\displaystyle
\varphi_{L_1}(t)=\int_{-\infty}^{\infty}\mathrm{e}^{\mathrm{i}tx}
g_1(x)\,\mathrm{d}x= \int_{-\infty}^{\ln A}\mathrm{e}^{\mathrm{i}tx}
\frac{\mathrm{e}^x}{A}\,\mathrm{d}x=\frac{\mathrm{e}^{\mathrm{i}t\,\ln
A}}{1+\mathrm{i}t}$ und der Produktformel folgt
$\displaystyle \varphi_{L_n}(t)=\left(\frac{\mathrm{e}^{\mathrm{i}t\,\ln
A}}{1+\mathrm{i}t}\right)^n$ und daher f\"ur $\displaystyle n\geq 2$ (als
\textit{uneigentliches} Integral auch f\"ur $\displaystyle n=1$)
\[ g_n(x)=\frac{1}{2\pi}\int_{-\infty}^{\infty} 
\mathrm{e}^{-\mathrm{i}tx} 
\left(\frac{\mathrm{e}^{\mathrm{i}t\,\ln A}}
{1+\mathrm{i}t}\right)^n\,\mathrm{d}t=
\frac{1}{2\pi}\int_{-\infty}^{\infty} 
\frac{\mathrm{e}^{\mathrm{i}t(n \ln A-x)}}
{(1+\mathrm{i}t)^n}\,\mathrm{d}t.
\]
Mit dem Residuensatz erhalten wir f\"ur $x\leq n \ln A$
\[
g_n(x)\!=\!\frac{1}{2\pi}\cdot2\pi\mathrm{i}\,
\underset{z=\mathrm{i}}{\mathrm{Res}}\frac{\mathrm{e}^
{\mathrm{i}z(n \ln A-x)}}{(1+\mathrm{i}z)^n}
\!=
\!\frac{\mathrm{i}}{(n\!-\!1)!}\left.\frac{\mathrm{d}^{n-1}}{\mathrm{d}z^{n-1}}
 \!\! \left(\frac{(z\!-\!\mathrm{i})^n}{(1\!+\!\mathrm{i}z)^n}\,
\mathrm{e}^{\mathrm{i}z(n \ln A -x)}\right)\right|_{z=\mathrm{i}}\!\!\!
=\!\frac{\mathrm{e}^x(n \ln A \!-\!x)^{n-1}}{A^n\,(n-1)!}
\]
Da im Falle \ $x>n\,\ln A$ \ \"uber einen \ Halbkreisbogen in
der \textit{unteren} Halbebene zu inte\-grieren ist (Lemma von Jordan),
folgt \ \ $g_n(x)=0$ \ f\"ur diese $x$.

Die Varianz von $\ln Y_n$: 
\[ \mathrm{Var}(\ln Y_1)=
\int_{-\infty}^{\ln A} x^2\,g_1(x)\,\mathrm{d}x-\bigl(\mathrm{E}(\ln
Y_1)\bigr)^2
=\int_{-\infty}^{\ln A} x^2\frac{\mathrm{e}^x}{A}\,\mathrm{d}x-(\ln
A-1)^2=1,
\]
folglich  \ \ $ \boldsymbol{\mathrm{Var}(\ln Y_n)= n}$
\ \ wegen Unabh\"angigkeit der Summanden $\displaystyle \ln X_i$.

Wir analysieren die auf \ \ $ \mu=0,  \ \ \sigma=1$
\ \textit{umskalierten} (standardisierten) \ Gr\"o{\ss}en\\
$\displaystyle \boldsymbol{\tilde{L}_n :=\frac{\ln Y_n-n(\ln
A-1)}{\sqrt{n}}}$ und ihre Dichten \ $\displaystyle
\tilde{g}_n(x)=\frac{1}{a}g_n\bigl(\frac{x-b}{a}\bigr)$ mit $\displaystyle
a=\frac{1}{\sqrt{n}}$ und $\displaystyle b= -\sqrt{n}(\ln A-1)$. \ \ Also folgt
\[ \tilde{g}_n(x)=\sqrt{n}\,g_n\bigl(x\sqrt{n}\!+\!n(\ln
A\!-\!1)\bigr)
=\sqrt{n}\frac{\mathrm{e}^{x\sqrt{n}+n(\ln A-1)}
(n\!-\!x\sqrt{n})^{n-1}}{A^n\,(n-1)!}=\frac{\sqrt{n}\,n^n}{n!\,\mathrm{e}^n}
\mathrm{e}^{x \sqrt{n}}\left(1\!-\!\frac{x \sqrt{n}}{n}\right)^{n-1}.
\]
Mit der Formel von \textit{de Moivre} und \textit{Stirling}, also
$\displaystyle n!=\sqrt{2\pi
n}\Bigl(\frac{n}{\mathrm{e}}\Bigr)^n\mathrm{e}^{1/(12n+\varepsilon_n)}$
mit$\displaystyle \varepsilon_n\to 0$ f"ur 
 $n\to \infty$ \ (und $
0<\varepsilon_n<1 \ (n\in\mathbb{N}), \ \ \varepsilon_n<0.1 \ (n\geq 4)$),
ergibt sich \\
$\displaystyle
\tilde{g}_n(x)=\frac{1}{\sqrt{2\pi}}\,\mathrm{e}^{-1/(12n+\varepsilon_n)+x\sqrt{n}+(n-1)\ln(1-x/\sqrt{n})}$.
\ \  Es gilt 
\[
(n\!-\!1)\ln(1\!-\!x/\sqrt{n})=-(n\!-\!1)\bigl(\frac{x}{\sqrt{n}}\!+\!\frac{x^2}{2n}\!+
\!\frac{x^3}{3n\sqrt{n}}\!+\cdots\bigr)=
-\!x\sqrt{n}\!+\!\frac{x}{\sqrt{n}}\!-\!\frac{x^2}{2}\!+\!\frac{x^2}{2n}\!+\!\frac{(1\!-\!n)x^3}{3n\sqrt{n}}\!+\cdots
\]
und damit f\"ur $|x|\leq c$ mit festem $\displaystyle c>0$
\boldmath\[
\tilde{g}_n(x)=\frac{1}{\sqrt{2\pi}}\,
\mathrm{e}^{-x^2/2+\mathcal{O}(1/\sqrt{n})} \ (n\to\infty).
\]\unboldmath%
Ein \textit{leicht} durchzurechnender Spezialfall des \textit{zentralen
Grenzwertsatzes}, analog zum klassi\-schen Satz von de Moivre und Laplace,
bei dem die Umfor\-mungen aber etwas verwickelter sind. (Genau um diesen
Satz zu beweisen, stell\-te ja auch Abra\-ham de Moivre die
"'Stirlingsche Formel"`
auf, bei der Stir\-ling nur den
kon\-stan\-ten Faktor mittels des Wallis-Produktes bestimmte.)

Die anfangs formulierte Ad-hoc-Erkl\"arung des den Studenten verbl\"uffenden
Ph\"a\-nomens ist damit analytisch untermauert.

\bigskip 

\textbf{IV. d) \ } N\"ahrungsweise Berechnung des Median

Ausgangspunkt war die Diskrepanz zwischen Erwartungswert und \textit{Median}
bei den Produkten. Daher geht es nun darum, letzteren analytisch
darzustellen. Was sich als nicht ganz einfach erweisen wird.

Wir diskutieren zun\"achst die beiden einfachsten F\"alle $n=2$
und $n=3$.

Der Median des Produktes zweier Zufallszahlen aus $\left(0,
A\right]$ ist charakterisiert als der Wert $c>0$, f\"ur den
die Hyperbel $x_1 x_2=c$ das Quadrat $0\leq x_1,
x_2 \leq A$ hal\-biert. 
Es ergibt sich $A^2/2=c(1+2\ln A -\ln c)$, also mit
$\tilde{c}:=c/A^2$: \ \  $
\boldsymbol{1/2=\tilde{c}(1-\ln \tilde{c})}$.
Dabei $0<\tilde{c}<1$, da ja offenbar $0<c<A^2$;
genauer folgt \ $\tilde{c}<1/2$.
Folglich $\displaystyle \tilde{c}=\frac{-1/2}{W_u(-1/2\mathrm{e})}$, wobei
$W_u(x)$ den \textit{unteren} Zweig der \textit{Lambertschen
W-Funktion} bezeichne, der Umkehrung von $
f(x)=x\,\mathrm{e}^x$, definiert f\"ur $ x\geq -1/\mathrm{e}$.
Der untere Zweig existiert nur f\"ur $ -1/\mathrm{e}\!\leq\!
x\!<\!0$, der obere Zweig $ W(x)$ f\"ur alle $
x\!\geq\! -1/\mathrm{e}$. Es ergibt sich \ \ $\displaystyle
W_u(-1/2\mathrm{e})=-2.678346990\ldots$ und \ \ $\displaystyle
\boldsymbol{\tilde{c}=0.18668230\ldots}$

Ein Plot, bei dem die kr\"aftigere Hyperbel den Median kennzeichnet und das
Quadrat halbiert, w\"ahrend \ die schw\"achere den Erwartungswert
repr\"asentiert:

{\centering   
\includegraphics[width=7.0cm,height=7.0cm]{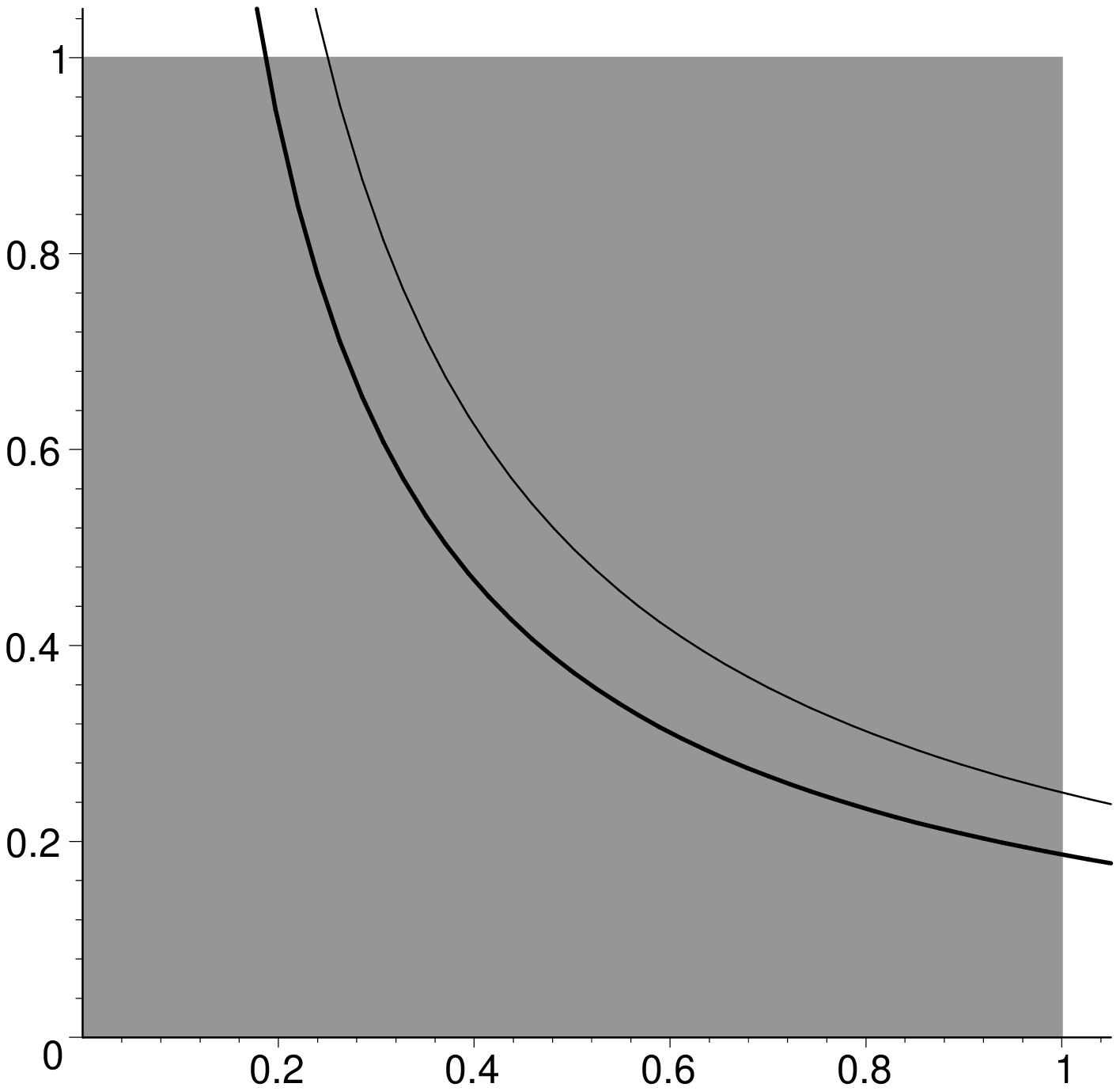}
\par}

Analog der Median bei drei Faktoren: Wir bestimmen den Wert $
c$, f\"ur welchen die Fl\"ache $ x_1 x_2 x_3=c$ das Volumen des
W\"urfels $ 0\leq x_1, x_2, x_3 \leq A$ in zwei
gleich\-gro{\ss}e H\"alften zerlegt.

Das ergibt \ $ \boldsymbol{c\approx0.0689716 A^3}$, \ wobei
$ \tilde{c}:=c/A^3$ der Gleichung 
$ \boldsymbol{1/2=\tilde{c}(1-\ln \tilde{c}+(\ln \tilde{c})^2/2)}$ \ gen\"ugt. 
\ 
Die grafische Veranschaulichung:

\includegraphics[width=8.2cm,height=9.0cm]{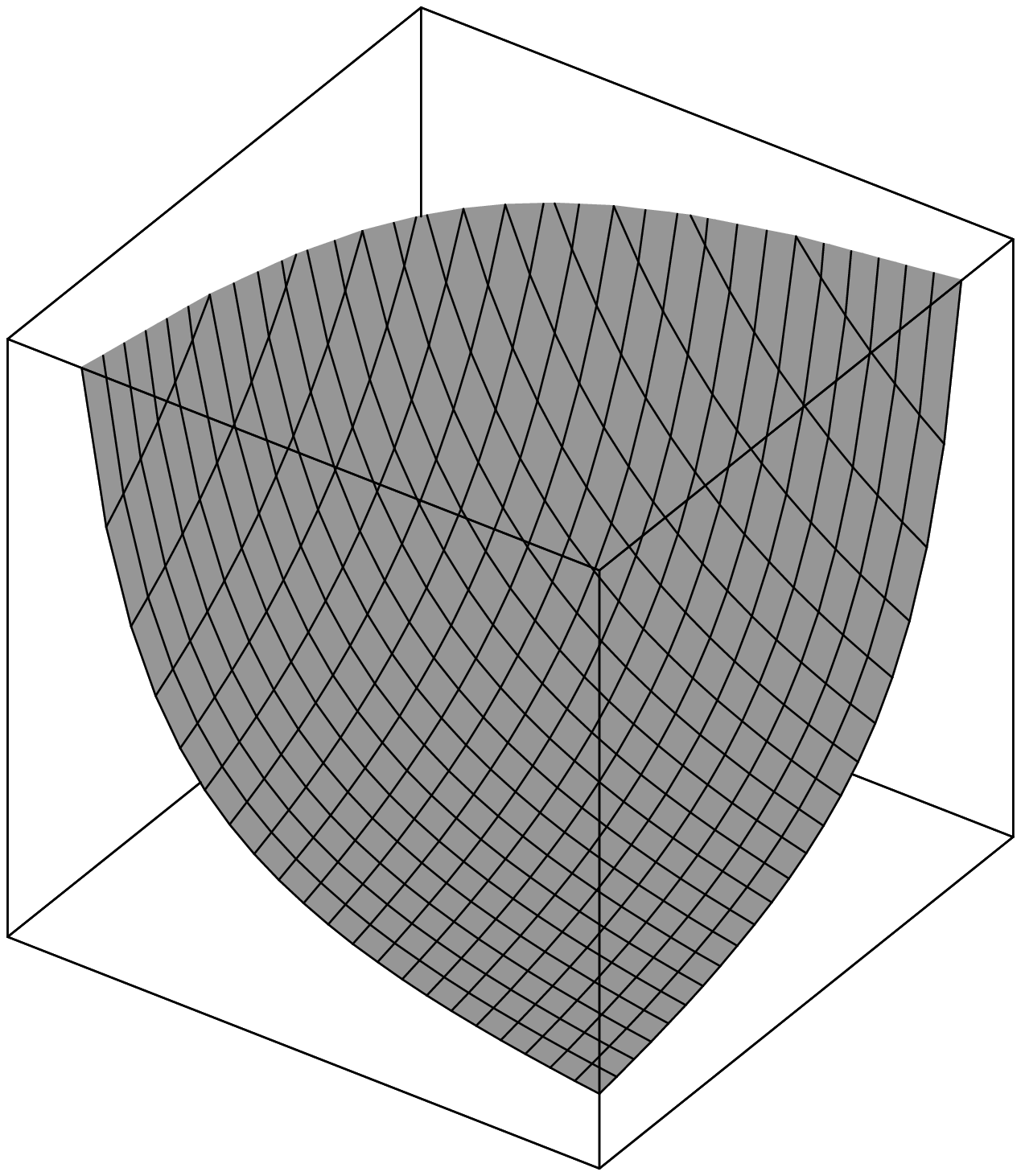}
\ \
\includegraphics[width=8.2cm,height=9.0cm]{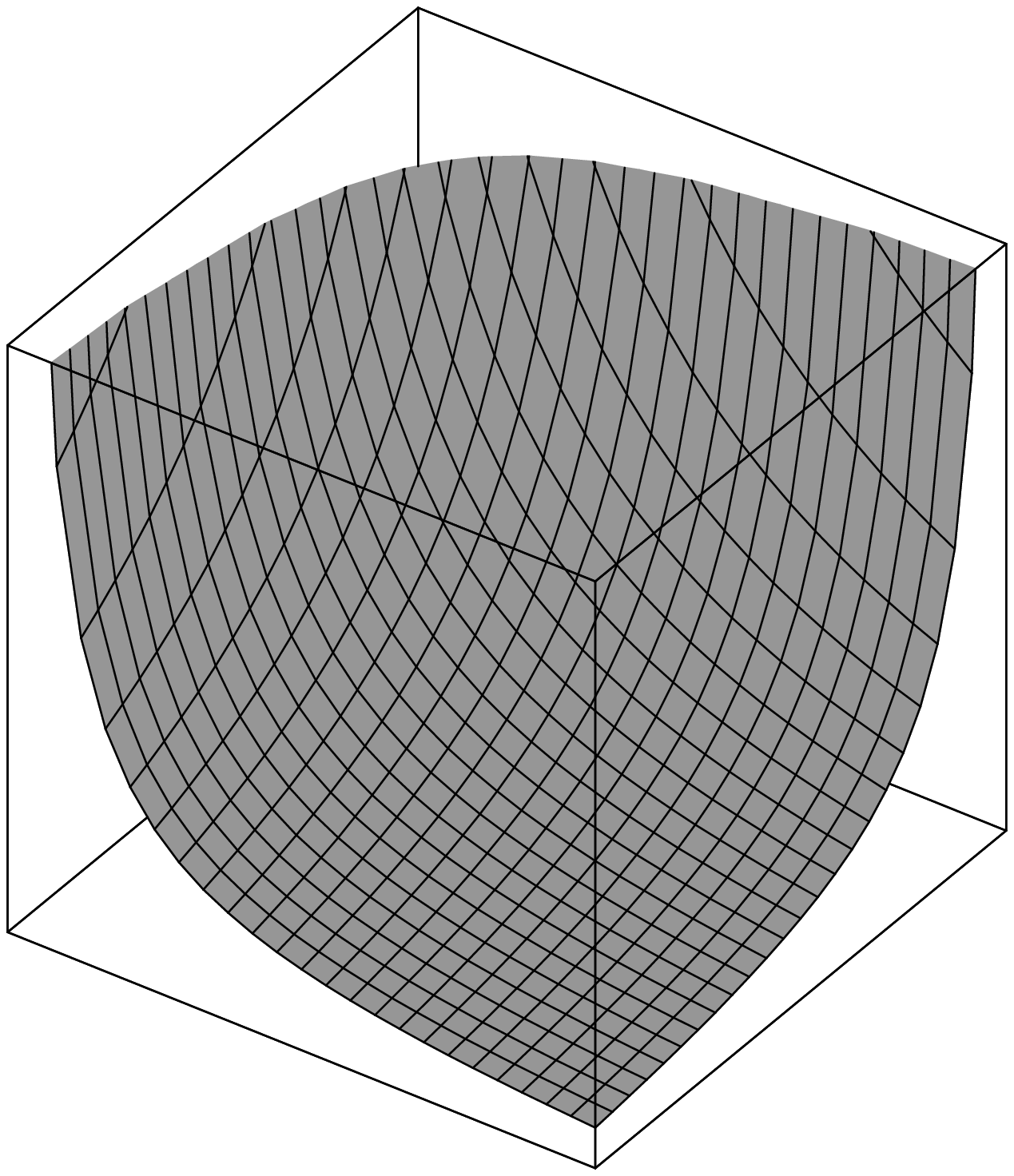}

Die Erwartungswert-Fl\"ache $ x_1 x_2 x_3 = A^3/8$ (linke
Grafik) trennt nur ein gutes Drittel des W\"urfel-Volumens ab.

Nun zum Allgemeinfall. Zun\"achst gilt 
\[ \mathrm{P}\bigl(\ln Y_n \leq n(\ln
A-1)\bigr)=\int_{-\infty}^0\tilde{g}_n(x)\,\mathrm{d}x
=\int_{-\infty}^0\frac{\sqrt{n}\,n^n}{n!\,\mathrm{e}^n}
\mathrm{e}^{x \sqrt{n}}\left(1\!-\!\frac{x \sqrt{n}}{n}\right)^{n-1}
\mathrm{d}x.
\]
Indem wir die Reihenentwicklung des Logarithmus und die Stirling\-sche Reihe
heran\-ziehen, erhalten wir nach \textit{einigem}
Neben\-rech\-nungs\-auf\-wand
\boldmath\[
\mathrm{P}(Y_n\leq (A/\mathrm{e})^n)
=\frac{1}{2}-\frac{1}{\sqrt{2\pi n}}\left(\frac{1}{3}+
\frac{1}{540n}-\frac{25}{6048 n^2} -
\frac{101}{155520 n^3}\right)
%+\frac{3184811}{3695155200n^4}+
%\frac{2745493}{8151736320 n^5}\right)+O\left(\frac{1}{n^6}\right)
+\mathcal{O}\left(\frac{1}{n^{9/2}}\right).
\]\unboldmath%
Andererseits 
\[\begin{split}
&\mathrm{P}\bigl( (A/\mathrm{e})^n<X_n\leq  (A/\mathrm{e})^n
(1+\varepsilon)\bigr)=
\int_{(A/\mathrm{e})^n}^{(A/\mathrm{e})^n(1+\varepsilon)} 
\frac{(n \ln A - \ln
  x)^{n-1}}{(n-1)!\cdot A^n} \,\mathrm{d}x
\\
&\phantom{P((A}{}=\int_1^{1+\varepsilon}
\frac{(n-\ln t)^{n-1}}{(n-1)!\cdot \mathrm{e}^n} \,\mathrm{d}t
 =\frac{(n/\mathrm{e})^n}{n!}\cdot\int_1^{1+\varepsilon}
\left(1-\frac{\ln t}{n}
\right)^{n-1}\,\mathrm{d}t.
\end{split}\]
Nun ist \ $ \varepsilon$ \ so zu w\"ahlen, dass $\displaystyle
\int_1^{1+\varepsilon}\left(1-\frac{\ln t}{n}\right)^{n-1}\mathrm{d}t$
m\"oglichst genau mit 
{\small\[
\frac{n!}{\sqrt{2\pi n}\cdot(n/\mathrm{e})^n}\left(\frac{1}{3}+
\frac{1}{540n}-\frac{25}{6048 n^2} -
\frac{101}{155520 n^3}\right)+\mathcal{O}\left(\frac{1}{n^4}\right)=
\frac{1}{3}+\frac{4}{135n}-\frac{8}{2835n^2}-\frac{16}{8505n^3}+\cdots
\]}
\"ubereinstimmt. 

So folgt $ \ln(1\!+\!\varepsilon)\!=\!1/3+\ln(1\!+\!a/n)$,
$ a=-8/405\!+\!b/n$, $ b=-8056/1148175+c/n$, usw.

Damit ist der Median $ M_n$ des Produkts $ Y_n$ von
$ n$ unabh\"angigen auf $ \left(0,A\right]$
gleichverteilten Zufallszahlen gegeben durch
\boldmath\[
\boxed{M_n=\sqrt[3]{\mathrm{e}}\;\left(\frac{A}{\mathrm{e}}\right)^n\cdot\left(1-\frac{8}{405\,
n}
-\frac{8056}{1148175\, n^2}+\mathcal{O}(n^{-3})\right)}
\]\unboldmath%
F\"alle $\displaystyle n\!=\!2$ und $\displaystyle n\!=\!3$: Es gilt
$\displaystyle
\sqrt[3]{\mathrm{e}}\Bigl(\frac{A}{\mathrm{e}}\Bigr)^2\!\!\approx 0.1888756
A^2$ und $\displaystyle
\sqrt[3]{\mathrm{e}}\Bigl(\frac{A}{\mathrm{e}}\Bigr)^3\!\!
\approx 0.06948345 A^3$. Allein der \textit{Hauptterm} stimmt also schon in
diesen einfachsten F\"allen recht gut mit den vorher angegebenen exakten
Median-Werten \"uberein.

\bigskip

\textbf{IV. e) \ } 50000 Produkte von je 500 Zufallszahlen

Ein gr\"o{\ss}eres Testbeispiel zu $\displaystyle Y_{500}$, generiert mit dem
Statistik-System \textbf{R}.

\verb!> f1<-function(A,n) prod(runif(n,0,A))}!
definiert eine Funktion, die $n$ gleichverteilte Zufallszahlen
aus dem Intervall $(0,A]$ erzeugt und dann deren Produkt
berechnet. 

Mit \verb!> f2<-function(A,n,r) {x<-numeric();for(k in 1:r) {x<-c(x,f1(A,n))};x}! 
und \verb!> Y<-f2(exp(1),500,50000)! erzeugt man einen
Datenvektor aus 50000 Pro\-duk\-ten von je\-weils 500 gleichverteilten
Zufallszahlen aus dem Intervall $(0, \mathrm{e}]$. Der R-Befehls "`summary"' liefert einen \"Uberblick:\\[-24pt]
\begin{verbatim}
> summary(Y)}
     Min.   1st Qu.    Median      Mean   3rd Qu.      Max.
1.668e-41 3.592e-07 1.395e+00 2.727e+33 4.917e+06 7.136e+37

> L<-log(Y)
> summary(L)
     Min.   1st Qu.    Median      Mean   3rd Qu.      Max.
-93.89000 -14.84000   0.33270   0.05403 15.41000   87.16000

>  hist(L,prob=TRUE);lines(density(L,bw=0.5),lwd=2)
\end{verbatim}
{\centering\includegraphics[width=12cm,height=8cm]{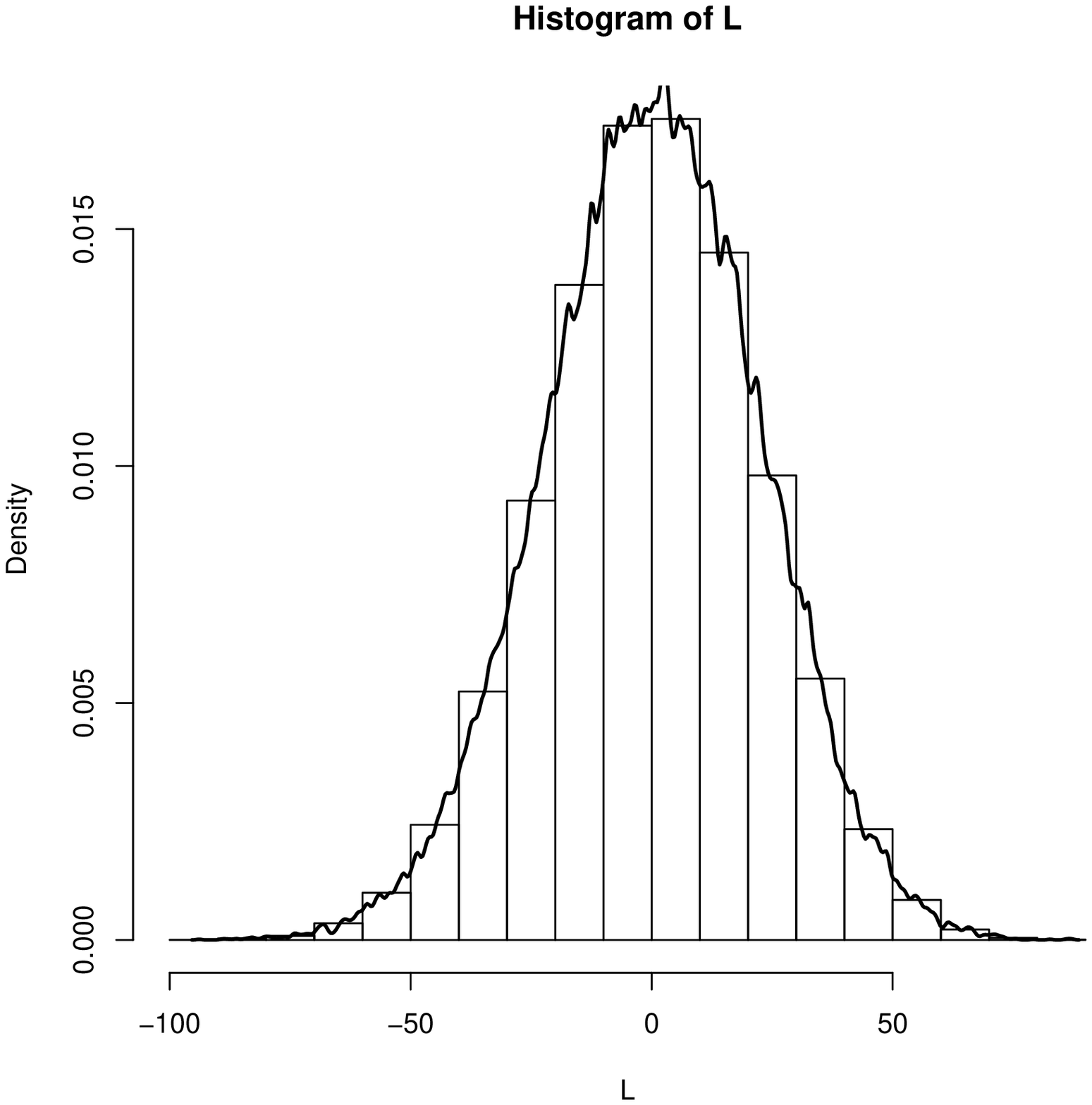}\par}

\end{document}